\documentclass[article]{siamart1116}[11pt]

\usepackage{amsmath,amssymb}
\usepackage{epsfig,graphics}
\usepackage{euscript,multicol, subfigure}
\usepackage{subfigure,sidecap}
\usepackage{wrapfig, comment}
\usepackage{colortbl,hhline}
\usepackage{mathtools}
\graphicspath{{./}}

\usepackage{pgfplots}
\usetikzlibrary{arrows,shapes}
\pgfplotsset{compat=1.3}

\usepackage{ dsfont }
\usepackage[textsize=footnotesize]{todonotes}

\renewcommand{\thesection}{\arabic{section}}

\usepackage{xcolor}

\DeclareMathAlphabet{\mathpzc}{OT1}{pzc}{m}{it} 

\usepackage{graphicx}
\usepackage{algorithm}
\usepackage{algorithmic}
\usepackage{pgf,tikz,pgfplots,import}
\usepackage{float}
\pgfplotsset{compat=newest}
\pgfplotsset{plot coordinates/math parser=false}
\newlength\figureheight
\newlength\figurewidth
\usetikzlibrary{plotmarks}
\usetikzlibrary{pgfplots.groupplots}
\usetikzlibrary{matrix}
\usepgflibrary{fadings}

\usepgfmodule{oo}

\makeatletter
\tikzset{
    every picture/.style={
        execute at begin picture={
            \let\ref\@refstar
        }
    }
}
\makeatother

\definecolor{plt-lc1}{rgb}{0.0078,0.2980,0.7961}
\definecolor{plt-lc2}{rgb}{1.0000,0.6431,0.2627}
\definecolor{plt-lc3}{rgb}{1.0000,0.2863,0.5255}
\definecolor{plt-lc4}{rgb}{0.6118,0.1765,1.0000}

\newlength\commLen
\newcommand{\mcaption}[3]{
  \ifthenelse{\isempty{#2}}
             {\caption{\textit{#3}\label{#1}}}
             {\caption[#2]{\textit{{\sc #2.}~#3}\label{#1}}}
             }

\newcommand{\algcaption}[3]{
        \ifthenelse{\isempty{#3}}
                   {\caption[#1]{{\sc #2.} \label{#1}}}
                   {\caption[#1]{{\sc #2.} \newline\small{#3} \label{#1}}}
        }

\makeatletter
\newcommand{\hypertargetraised}[1]{\Hy@raisedlink{\hypertarget{#1}{}}}
\makeatother

 \newcommand\authornote[2][Note]{}

\newcommand{\bbm}{\begin{bmatrix}}
\newcommand{\ebm}{\end{bmatrix}}

\newcommand\lbl[1]{\ensuremath{\mathrm{#1}}}

\newcommand\sci[2][1]{\ensuremath{\ifthenelse{\equal{#1}{1}}{}{#1\times}10^{#2}}}
\newcommand{\cA[1]}{\mathpzc{#1}}
\newcommand{\ttA[1]}{\mathtt{#1}}
\newcommand{\acc}{\varepsilon}
\newcommand{\vc[1]}{\boldsymbol{#1}}

\def\defeq{\mathrel{\mathop :}=}
\newcommand\tensor[1]{\mathpzc{#1}}

\newcommand\tree[1]{\mathcal{#1}}
\newcommand\la[1]{\mathsf{#1}} 
\newcommand\laT[1]{\la{#1}}    
\newcommand\ordinal[1]{\ensuremath{#1\mathrm{th}}}

\def\vec{\ensuremath{\mathrm{vec}}}
\def\reshape{\ensuremath{\mathtt{reshape}}}

\def\TT{TT}
\def\TTrain{tensor train}

\def\Mdash/{\kern 0.05em---\kern 0.05em}

\title{Tensor Train accelerated solvers for nonsmooth rigid body dynamics}
\author{Eduardo Corona~\footnotemark[1], \ David Gorsich~\footnotemark[2] , \ Paramsothy Jayakumar~\footnotemark[2]  \ and  Shravan Veerapaneni~\footnotemark[1]}
\date{}
\begin{document}
\maketitle

\renewcommand{\thefootnote}{\fnsymbol{footnote}}
\footnotetext[1]{Department of Mathematics, University of Michigan}
\footnotetext[2]{U.S. Army TARDEC}

\begin{abstract} 
In the last two decades, increased need for high-fidelity simulations of the time evolution and propagation of forces in granular media has spurred a renewed interest in the discrete element method (DEM) modeling of frictional contact. Force penalty methods, while economic and widely accessible, introduce artificial stiffness, requiring small time steps to retain numerical stability. Optimization-based methods, which enforce contacts geometrically through complementarity constraints leading to a differential variational inequality problem (DVI), allow for the use of larger time steps at the expense of solving a nonlinear complementarity problem (NCP) each time step. We review the latest efforts to produce solvers for this NCP, focusing on its relaxation to a cone complementarity problem (CCP) and solution via an equivalent quadratic optimization problem with conic constraints. We distinguish between \emph{first order} methods, which use only gradient information and are thus linearly convergent and \emph{second order} methods, which rely on a Newton type step to gain quadratic convergence and are typically more robust and problem-independent. However, they require the approximate solution of large sparse linear systems, thus losing their competitive advantages in large scale problems due to computational cost.

In this work, we propose a novel acceleration for the solution of Newton step linear systems in second order methods using low-rank compression based fast direct solvers, leveraging on recent direct solver techniques for structured linear systems arising from differential and integral equations. We employ the Quantized Tensor Train (QTT) decomposition to produce efficient approximate representations of the system matrix and its inverse. This provides a versatile and robust framework to accelerate its solution using this inverse in a direct or a preconditioned iterative method. We demonstrate compressibility of the Newton step matrices in Primal Dual Interior Point (PDIP) methods as applied to the multibody dynamics problem. Using a number of numerical tests, we demonstrate that this approach displays \emph{sublinear} scaling of precomputation costs, may be efficiently updated across Newton iterations as well as across simulation time steps, and leads to a fast, optimal complexity solution of the Newton step. This allows our method to gain an order of magnitude speedups over state-of-the-art preconditioning techniques for moderate to large-scale systems, hence mitigating the computational bottleneck of second order methods.
\end{abstract}

\section{Introduction}
\label{sec:intro}
The discrete element method (DEM) \cite{cundall1979DEM} is one of the most widely used approaches to simulate the multi-body dynamics such as in granular materials. This method considers the granular medium as a collection of discrete particles; each responding to body forces such as gravity, inertia or drag, as well as repulsive or dissipative forces caused by contact. Materials in granular form are omnipresent in industry and understanding their dynamics is crucial for a broad range of application fields, including terramechanics, active media, additive manufacturing, nanoparticle self-assembly, avalanche dynamics, composite materials, pyroclastic flows, etc. 

Distinct instances of the DEM are defined essentially by their modeling of contact forces, and thus in how collisions are handled. We distinguish two main types in the literature: force penalty methods (DEM-P) and complementarity formulations (DEM-C). Penalty methods introduce one or multiple layers of spring-like forces between objects in contact, and may introduce additional fields to represent friction. These methods are widely used due to being computationally inexpensive and easy to implement. However, in many cases they introduce high stiffness, requiring extremely small time steps to retain stability during collisions. An in-depth and informative comparison between the two can be found in  \cite{pazouki2017PversusC}. 

DEM-C methods such as in \cite{tasora2010DEMC,anitescu2006optimization} enforce contacts using complementarity constraints, leading to a differential variational inequality (DVI) problem upon discretization. This allows for the use of larger time steps in its integration, as contacts are enforced geometrically. Given a time-stepping scheme for this problem, a nonlinear complementarity problem (NCP) must be solved to compute the corresponding contact forces at each time step. One common way to solve this optimization problem is via a linearization of the constraints, producing a linear complementarity problem (LCP) \cite{stewart1996lcp,anitescu1997lcp}. An alternative relaxation method produces an equivalent, convex cone complementarity problem (CCP). Efficient quadratic cone programming techniques such as those in \cite{tasora2010DEMC,mazhar2015Nesterov,Negrut2013Krylov,petra2009computational,fang2015primal,kleinert2015simulating} have been proposed to solve the CCP. 

In \cite{melanz2017comparison}, the authors performed a comparison of these quadratic programming techniques focused on determining which class of methods performed best: second order optimization methods, i.e. those that use Hessian information (Interior Point), or first order methods, i.e. methods using only gradient information (Jacobi, Gauss-Seidel, Projected Gradient Descent). Second order methods display quadratic convergence in a neighborhood of the solution, and thus their convergence is much faster and more problem-independent than that of linearly convergent first order methods, requiring at least 1 to 2 orders of magnitude less iterations to reach a desired accuracy across experiments in this work. However, each iteration in second order metods requires the solution of a (generally sparse) linear system, which can become costly as the dimensions of the many-body problem increase. For large problems, the added computational effort ultimately eclipses the gains obtained from the reduction in iteration counts.     

The preferred solution method for large sparse systems is often a Krylov subspace iterative method. The iteration counts, and thus the performance of these methods are known to be directly affected by the eigenspectrum of the associated matrices. Preconditioning techniques can be used to cluster the eigenvalues away from zero and drastically reduce iteration counts; we refer the reader to \cite{benzi2002preconditioning} for a general review. Generic sparse preconditioners are most typically based on incomplete or sparsified factorizations, such as the well-known Incomplete LU and Cholesky methods \cite{saad2003iterative}. In \cite{melanz2017comparison}, a fast, parallel SaP (split and paralellize) \cite{SaPGPU} preconditioner was used to accelerate the PDIP method, garnering reductions in iteration counts and execution times.  

Most general-purpose preconditioners suffer a trade-off between precomputation costs and the resulting reduction in iteration counts, and their performance is often problem-dependent. Moreover, in the context of optimization problems such as the ones we are interested in, system matrices change every iteration, and it is often impossible or expensive to update the associated preconditioners. Finally, improving their scaling with the number of degrees of freedom and the level of sparsity is extremely important, as it is central to remaining competitive in large-scale problems. 

In the last decades, a continued effort has been made to produce \textit{direct} solvers for structured linear systems arising from differential and integral equations. These solvers entirely side-step the challenges related to convergence speed of iterative solvers. They can also lead to dramatic improvements in speed, in particular in situations where a large number of linear systems with coefficient matrices that stay fixed or can be updated via low rank modifications. Additionally, they also provide a methodology to produce robust preconditioners: low accuracy direct solvers may be used in conjunction with iterative refinement or the Krylov subspace method of choice; the use of an approximate low accuracy inverse requires less memory and precomputation time than a direct solver, at the expense of a slight increase in the number of iterations. These features have led to the adoption of these solvers in other areas of scientific computing and statistics.

In \cite{corona2015tensor}, an effective and memory-efficient solver based on the quantized tensor train decomposition (QTT) was presented. By recasting system matrices as tensors, the tensor train TT compression and inversion routines were used to produce direct solvers and robust preconditioners for integral equations in complex geometries in three dimensions. Key properties of this solver that differentiate it from other hierarchical matrix approaches feature sublinear computational costs and memory requirements with problem size $N$, as well as techniques to produce economic updates for a matrix and its inverse across time-steps, \emph{even when matrix size changes}. 

The main goal of this work is to employ a QTT-based approach to provide a radical speed up to second order optimization methods. We demonstrate its application to interior point methods such as the PDIP in the context of many-body dynamics problems; however, we expect our discussion to apply with little to no modification to a general class of Newton and Quasi-Newton type methods. We first study compressibilty of the system matrices and their inverses in TT format for a range of target accuracies. We then demonstrate how factorization re-use can provide significant speed-ups to precomputation costs, reducing costs by orders of magnitude. For three validation tests of common soil mechanics phenomena---sedimentation, blade drafting and direct shear experiments---we show that a TT-based preconditioner displays efficient and robust performance for problems with $>10^4$ bodies, garnering up to an order of magnitude speed-up and greatly improved iteration counts when compared with state-of-the-art ILU-preconditioned methods. We also confirm an extremely significant gain in scaling of precomputation costs: precomputation for the TT preconditioner is \emph{sublinear} (for all practical purposes, constant) as the number of collisions and matrix size increase.  

The QTT decomposition is one of several approaches for approximate solution of linear systems based in hierarchical compression. In the context of sparse structured systems and related factorizations, work has been done for a number of hierarchical matrix formats: \textit{Hierarchically Semi Separable (HSS)} \cite{G2011thesis,ambikasaran2013sparseHSS,xia2009multifrontal, Gu06sparse, Gu06ULV, xia2010,ho2015hierarchical_de}, \textit{$\mathcal{H}$ matrices:} \cite{borm2003hierarchical,   2008_bebendorf_book, 2010_borm_book}, and FMM \cite{coulier2017IFMM, pouransari2017fast}. Producing efficient, global factorization updates and dealing with high storage costs is an ongoing challenge in these alternate formats.

\section{Mathematical Preliminaries: dynamics of rigid bodies in the presence of friction}
\label{sec:formulation}
\subsection{Problem formulation\label{ssc:formulation}}

We consider a granular material comprised of $M$ rigid particles $B_i$ in $\mathbb{R}^3$; the position of each body is uniquely described by the coordinates $\vc[x]_i$ for its center of mass and the rotation $\vc[Q]_i$ of a reference frame fixed to the body, represented by a unimodular quaternion requiring three extra parameters. Let $\vc[q]_i = (\vc[x]_i,\vc[Q]_i)$ then be the $6$ generalized coordinates for the $i$-th body, and $\vc[q] = (\vc[q]_1,\vc[q]_2,\dots,\vc[q]_M) \in \mathbb{R}^{6M}$. 

Each particle's motion can thus be understood as a translation of its center, with velocity $\vc[u]_i$ and a rotation of its frame, with angular velocity $\vc[\omega]_i$. Applying Newton's second law, we relate the corresponding accelerations to the total force $\vc[F]_i$ and torque $\vc[T]_i$ applied to it. The general equations of motion are then given by 

\begin{align}
\vc[\dot{q}] &= \cA[L](\vc[q]) \vc[v] \label{eq:particle_evolution} \\
M(\vc[q]) \vc[\dot{v}] &= \vc[f]_B(\vc[q],\vc[v]) + \vc[f]_{C}. \label{eq:force_balance}
\end{align}

where $\cA[L](\vc[q])$ is a linear operator that relates velocities to the rate of change in generalized coordinates, $M$ is the mass matrix and $\vc[v], \vc[f]_B, \vc[f]_C \in \mathbb{R}^{6M}$ contain each body's translational and rotational velocities, the total body forces and torques applied to them and the reaction forces and torques due to contact dynamics, respectively. 

\subsection{Complementarity contact model\label{ssc:contact_ivp}}

We impose the following contact constraints: no two bodies should penetrate, and if there is contact, a normal force and a tangential frictional force act at the interface. Consider two bodies $B_{i_1}$ and $B_{i_2}$, and let $\Phi_i(\vc[q])$ be an unsigned distance function (also known as a gap function) for the pair $i  = (i_1,i_2)$ satisfying $\Phi_i(\vc[q])>0$ if the two bodies are separated, $\Phi_i(\vc[q])=0$ if they are touching, $\Phi_i(\vc[q])<0$ otherwise. 

If the pair of bodies touch ($\Phi_i(\vc[q])=0$), let $\vc[n]$,$\vc[t]_1,\vc[t]_2$ be unit normal and tangential vectors at the point of contact. Contact forces $f_N = \gamma_{i,n} \vc[n]$ and $f_T = \gamma_{i,1} \vc[t]_1 + \gamma_{i,2} \vc[t]_2$ are then applied to each body in opposite directions. The complementarity constraint for the normal force \eqref{eq:complementarity} prevents penetration, enforcing that bodies move away from each other at contact. The Coulomb friction model ties the magnitudes of the normal and tangential forces. Using a maximum dissipation principle, the friction force is posed as the solution to an optimization problem in $\eqref{eq:friction}$: 

\begin{align}
\gamma_{i,n} &\geq 0 \,\ \ \ \Phi_i(\vc[q]) \geq 0 \,\ \ \ \Phi_i(\vc[q]) \gamma_{i,n} = 0, \label{eq:complementarity} \\
(\gamma_{i,1},\gamma_{i,2}) &= \mathrm{argmin}_{||(\beta_1,\beta_2)|| \leq \mu \gamma_{i,n}} \vc[v]^T (\beta_{1} \vc[t]_1 + \beta_{2} \vc[t]_2), \label{eq:friction} 
\end{align}

where $\mu$ is the static friction coefficient. The feasible set of forces $f = f_N + f_T$ for the minimization problem in \eqref{eq:friction} are known as a \emph{friction cone} $\Upsilon_i = \{(\gamma_{i,n},\gamma_{i,1},\gamma_{i,2}) \ | \ ||(\gamma_{i,1},\gamma_{i,2})|| \leq \mu \gamma_{i,n} \}$. The complementarity condition in \eqref{eq:complementarity} is typically abbreviated using the notation $0 \leq \gamma_{i,n} \perp \Phi_i(\vc[q]) \geq 0$. We note that an alternate complementarity formulation for \eqref{eq:complementarity} can be obtained by replacing $\Phi_i$ with the normal velocity. 

Since all forces are zero in the absence of contact, we wish to incorporate contact constraints only for pairs of objects approaching collision (e.g. within the next timestep). For this purpose, we define the set $\cA[A]$ of pairs of bodies separated by a distance smaller than a threshold $\delta>0$ 

\begin{equation} \cA[A](\vc[q],\delta) = \{ i | \Phi_i(\vc[q]) \leq \delta \}. \end{equation}

For a system with $N_c = |\cA[A](q,\delta)|$ contacts, this adds a number of constraints to the equations of motion proportional to $N_c$. We note that for dense granular flows, $N_c$ itself scales as $O(M)$ for $M$ bodies. 

Now, let $\vc[D]_i = [D_{i,n} \ D_{i,1} \ D_{i,2}] \in \mathbb{R}^{6M \times 3}$ mapping the multipliers $\gamma_{i,n},\gamma_{i,1},\gamma_{i,2}$ to the forces and torques applied to bodies $B_{i_1}$ and $B_{i_2}$ at contact. Then, the equations of motion result in the  differential variational inequality

\begin{align}
\vc[\dot{q}] &= \cA[L](\vc[q]) \vc[v], \label{eq:particle_evolution_demc} \\
M(\vc[q]) \vc[\dot{v}] &= \vc[f]_B(\vc[q],\vc[v]) + \sum_{i \in \cA[A](\vc[q],\delta)} D_{i,n} \gamma_{i,n} + D_{i,1} \gamma_{i,1} + D_{i,2} \gamma_{i,2}, \label{eq:force_balance_demc} \\
i \in \cA[A](\vc[q],\delta) &: \gamma_{i,n} \geq 0 \ \perp \ \Phi_i(\vc[q]) \geq 0, \label{eq:complementarity_demc} \\
(\gamma_{i,1},\gamma_{i,2}) &= \mathrm{argmin}_{||(\beta_1,\beta_2)|| \leq \mu \gamma_{i,n}} \vc[v]^T (D_{i,1}\beta_{1} + D_{i,2}\beta_{2} ). \label{eq:friction_demc}  
\end{align} 

A semi-implicit, first order integration scheme is then used to advance this system in time. Given position $\vc[q]^k$ and velocity $\vc[v]^k$ at a given time step $t^k$ and step size $\Delta t$, velocity $\vc[v]^{k+1}$ and contact forces are solved via a nonlinear complementarity problem (NCP). The new velocity is used to evolve the position in time.  

\begin{align}
\vc[q]^{k+1} &= \vc[q]^k + \Delta t \cA[L](\vc[q]^k) \vc[v]^{k+1}, \label{eq:particle_evolution_disc} \\
M(\vc[q^{k}]) (\vc[v]^{k+1} - \vc[v]^{k}) &= \Delta t \vc[f]_B(\vc[q]^k,\vc[v]^k) + \sum_{i \in \cA[A](\vc[q]^k,\delta)} D_{i,n} \gamma_{i,n} + D_{i,1} \gamma_{i,1} + D_{i,2} \gamma_{i,2}, \label{eq:force_balance_disc} \\
i \in \cA[A](\vc[q]^k,\delta) &: \gamma_{i,n} \geq 0 \ \perp \ \frac{1}{\Delta t}\Phi_i(\vc[q]^k) + D^T_{i,n} \vc[v]^{k+1} \geq 0, \label{eq:complementarity_disc} \\
(\gamma_{i,1},\gamma_{i,2}) &= \mathrm{argmin}_{||(\beta_1,\beta_2)|| \leq \mu \gamma_{i,n}} (\vc[v]^{k+1})^T (D_{i,1}\beta_{1} + D_{i,2}\beta_{2} ). \label{eq:friction_disc}
\end{align} 

We note that \eqref{eq:complementarity_disc} is obtained from \eqref{eq:complementarity_demc} via a linearization, dividing by $\Delta t$ (which does not affect complementarity, but is numerically desirable). This makes the future velocity $v^{k+1}$ the sole variable needed to enforce the complementarity condition. We also note that in this discretization, $(\gamma_{i,n},\gamma_{i,1},\gamma_{i,2})$ constitute contact \emph{impulses}, i.e. force magnitudes multiplied by the step length $\Delta t$. 

\subsection{Solving the optimization problem\label{ssc:optimization}}

A relaxation over the complementarity constraint \eqref{eq:complementarity_disc} can be introduced \cite{anitescu2006optimization}, turning the problem into a convex, second-order cone complementarity problem (CCP). 

\begin{equation}
\gamma_{i,n} \geq 0 \ \perp \ \frac{1}{\Delta t}\Phi_i(\vc[q]^k) + D^T_{i,n} \vc[v]^{k+1} - \mu_i \sqrt{(D^T_{i,1}\vc[v]^{k+1})^2 + (D^T_{i,2}\vc[v]^{k+1})^2} \geq 0.
\label{eq:complementarity_relax}
\end{equation}

The solution of this relaxed problem approaches the solution of the NCP as the step-size $\Delta t$ goes to zero. Additionally, the CCP is equivalent to the KKT first-order optimality conditions for a quadratic optimization problem with conic constraints. We define a contact transformation matrix $\vc[D] = [D_1 \ ; D_2 \ ; \dots \ D_{N_c}] \in \mathbb{R}^{6M \times 3 N_c}$, a matrix $\vc[N]$ and vector $\vc[r]$: 

\begin{align}
\vc[N] &= \vc[D]^T \vc[M]^{-1} \vc[D], \label{eq:Ndef} \\
\vc[r]  &= \vc[b] + \vc[D]^T \vc[M]^{-1} \vc[k], \label{eq:rdef}
\end{align}

where $\vc[b] = [\vc[b]_1 \ ; \vc[b]_2 \ ; \dots \ \vc[b]_{N_c}]$, $\vc[b]_i = [\Phi_i^k / \Delta t \ ; 0 \ 0] \in \mathbb{R}^3$ and $\vc[k] = \vc[M] \vc[v]^k + \Delta t \vc[f]^k$. We note that $\vc[N]$ is a $3 N_c \times 3 N_c$ symmetric positive semi-definite matrix and typically sparse. The aforementioned quadratic program is then given by:

\begin{align}
\min q(\gamma) &= \frac{1}{2} \gamma^T \vc[N] \gamma + \vc[r]^T \gamma \label{eq:qp} \\
\mathrm{subject} \ \mathrm{to} \ \ \ & ||(\gamma_{i,1}, \gamma_{i,2})|| \leq \mu_i \gamma_{i,n} \ \ \ i =1,2,\dots,N_c. \label{eq:qp_constraints}
\end{align}  

While this relaxed problem may introduce artifacts when step size $\Delta t$, sliding velocity or friction are large, it enables the use of a wide range of quadratic programming methods for its solution: 

\begin{itemize}
\item Projected Jacobi and Gauss-Seidel methods \cite{tasora2010DEMC}.
\item Projected gradient descent methods like Accelerated Projected Gradient Descent \cite{mazhar2015Nesterov}, Barzilai - Borwein \cite{bbpgd} and the Kucera and Preconditioned spectral projected gradient with fallback (P-SPG-FB) methods in \cite{Negrut2013Krylov}. 
\item Krylov subspace methods: Gradient projected minimum residual (GPMINRES) in \cite{Negrut2013Krylov}. 
\item Primal-Dual Interior Point (PDIP) methods \cite{andersen2011interior,petra2009computational,fang2015primal}.
\item Symmetric Cone Interior Point (SCIP) methods \cite{kleinert2015simulating}.
\end{itemize}

Projected Jacobi and Gauss-Seidel methods, while requiring only fast matrix applies of $\vc[N]$, have slow, linear convergence, often requiring thousands of iterations to obtain a significant reduction for the objective function. The projected gradient descent and Krylov subspace methods have since been proposed, providing considerable iteration count reductions while retaining cost-efficiency per time step. However, they remain linearly convergent, and they require an increasing number of iterations as problem size increases for a variety of problems of interest. 

Interior point methods are often based on a modified Newton or Quasi-Newton step, displaying quadratic convergence near minima and iteration counts which are less dependent on problem size. However, they require the approximate solution of large linear systems in order to produce the Newton step, thus losing this competitive advantage due to computational cost.  

\subsection{Overview of interior point methods}

Interior point methods, also known as barrier methods, are a class of algorithms tailored to solve constrained convex optimization problems \cite{nocedal2006nonlinear}. That is,  

\begin{align}
\min & \ f_0 (\vc[\gamma]) \label{eq:cvx_objective}  \\
\mathrm{subject \ to} & \ f_i(\vc[\gamma]) \leq 0, \ i=1,\dots,m \label{eq:cvx_constraints}
\end{align}

for $f_i \in C^2(\mathbb{R}^n)$ and convex. They proceed by transforming this problem into an unconstrained minimization problem, encoding the feasible set defined by the constraints using a barrier function, e.g., the logarithmic barrier $B_t ( z ) = -(1/t) \log ( -z)$. They then pose the unconstrained convex problem: 

\begin{equation}
\min \ f_0 (\vc[\gamma]) + \sum_{i=1}^m B_t ( f_i(\vc[\gamma]) ). \label{eq:barrier_problem}
\end{equation}

The parameter $t$ controls the strength of the barrier, and as $t \rightarrow \infty$, $B(z) \rightarrow I(z)$ with $I(z)$ the indicator function over $(-\infty,0]$, and the problem in Eq \ref{eq:barrier_problem} becomes equivalent to the original constrained program defined by Eqs \ref{eq:cvx_objective} and \ref{eq:cvx_constraints}. Given the set of problems defined by Eq \ref{eq:barrier_problem}, interior point methods proceed by following along the corresponding \emph{central path} $\{ \vc[\gamma]^*(t) : t > 0 \}$ of optima, which are located inside of the feasible set, and converge to the solution of the original problem as $t \rightarrow \infty$.  

\subsubsection*{Primal-Dual Interior Point methods}

Primal-dual methods proceed by defining a path of solutions $(\vc[\gamma]^*, \vc[\lambda]^*)$ . These can be obtained by considering the KKT conditions of Eq \ref{eq:barrier_problem}, 

\begin{align}
f_i(\vc[\gamma]) &< 0, \ \ i=1,\dots,m \\
\nabla f_0(\vc[\gamma]) + \sum_{i=1}^m \frac{-1}{t f_i(\vc[\gamma])} \nabla f_i (\vc[\gamma] ) &= 0 
\end{align}

We then define the lagrange multipliers $\vc[\lambda]_i = \frac{-1}{t f_i(\vc[\gamma])}$, setting up the following conditions: 

\begin{align}
f_i(\vc[\gamma]) &< 0, \ \ i=1,\dots,m \\
\vc[\lambda]_i &< 0, \ \ i=1,\dots,m \\
-\vc[\lambda]_i f_i(\vc[\gamma]) &= \frac{1}{t}, \ \ i=1,\dots,m \label{eq:KKT_complementarity} \\ 
\nabla f_0(\vc[\gamma]) + \sum_{i=1}^m \vc[\lambda]_i \nabla f_i (\vc[\gamma] ) &= 0 \label{eq:KKT_lagrangian}. 
\end{align}

We note that, as $t \rightarrow \infty$, Eq \ref{eq:KKT_complementarity} becomes a strict complementarity condition. 

\subsubsection*{The PDIP step}

The PDIP step is then obtained by applying Newton's algorithm to solve these KKT conditions. Defining the residual of the system given by Eqs \eqref{eq:KKT_complementarity} and \eqref{eq:KKT_lagrangian} $\vc[r]_t (\vc[\gamma],\vc[\lambda])$ as: 

\begin{equation}
\vc[r]_t (\vc[\gamma],\vc[\lambda]) = \begin{bmatrix} \nabla f_0(\vc[\gamma]) + \nabla \vc[f](\vc[\gamma])^T \vc[\lambda] \\
-diag(\vc[\lambda]) \vc[f](\vc[\gamma]) - \frac{1}{t} \vc[1] \end{bmatrix} = 0
\end{equation}

the Newton step is then defined by the solution of the linear system given by: 

\begin{equation}
\begin{bmatrix} \nabla^2 f_0(\vc[\gamma]) + \sum_{i=1}^m \vc[\lambda]_i \nabla^2 \vc[f]_i (\vc[\gamma] ) & \nabla \vc[f] (\vc[\gamma] )^T \\
-diag(\vc[\lambda]) \nabla \vc[f](\vc[\gamma]) & -diag(\vc[f](\vc[\gamma]))  \end{bmatrix} 
\begin{bmatrix} \Delta \vc[\gamma] \\ \Delta \vc[\lambda] \end{bmatrix} = -\vc[r]_t (\vc[\gamma],\vc[\lambda]) \label{eq:cvx_Newton_system}
\end{equation}

this is typically coupled with a backtracking line search, and a strategy to increase $t$ until sufficient convergence to the solution of the original problem is satisfied. 

\subsubsection*{Application to the CCP} 
\begin{figure}
\begin{center}
\includegraphics[scale=0.25]{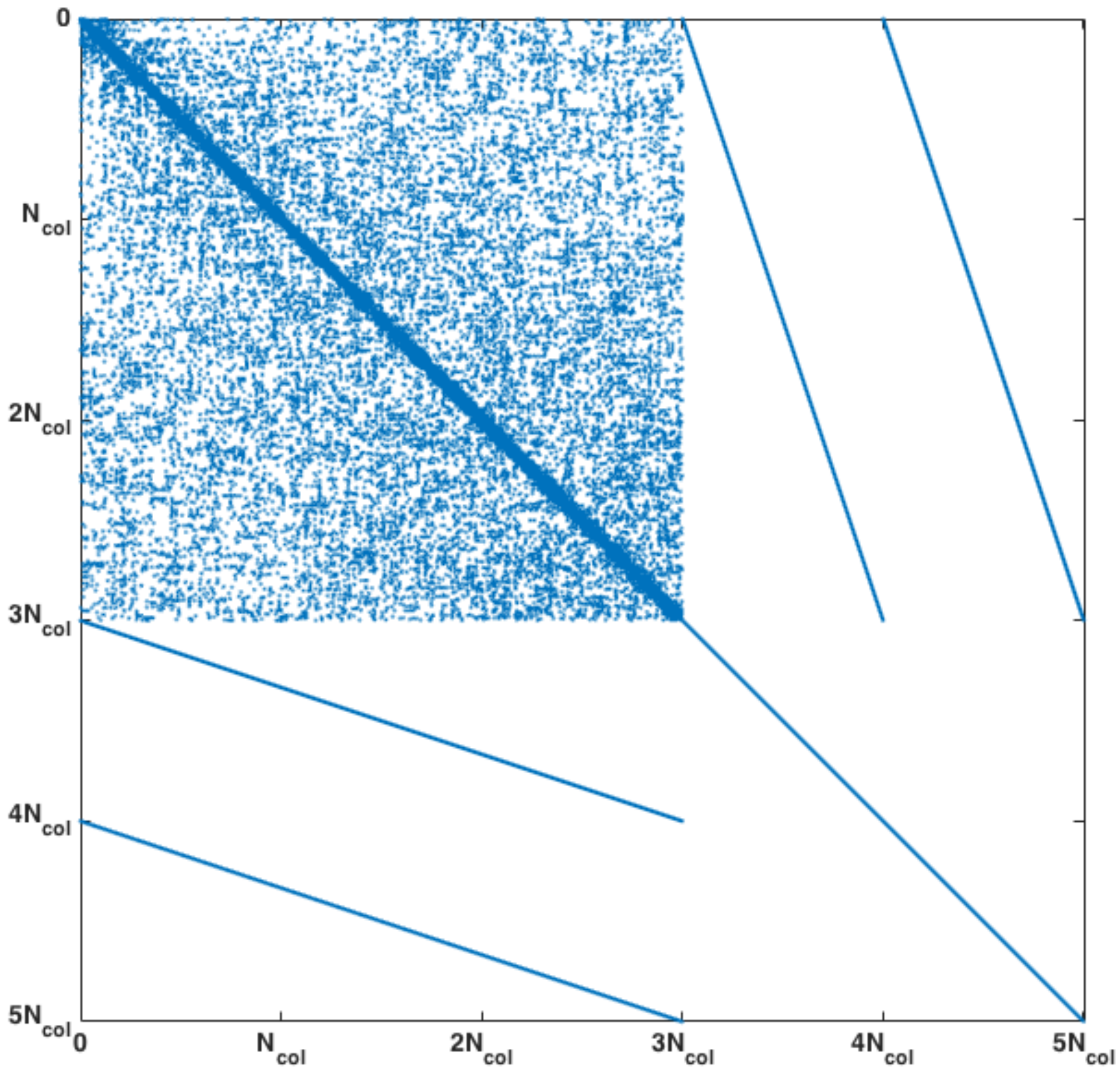}
\end{center}
\caption{Sparsity structure of PDIP Newton system matrix for a general multibody dynamics problem}
\label{fig:SpyNewtonMat}
\end{figure}

For the equivalent quadratic program with conic constraints in Eqs \eqref{eq:qp} - \eqref{eq:qp_constraints}, we define $f_0(\vc[\gamma]) = q(\vc[\gamma]) = \frac{1}{2} \gamma^T \vc[N] \gamma + \vc[r]^T \gamma$ and $2 N_c$ inequality constraints (Eq \ref{eq:qp_constraints}) given by $f_i(\vc[\gamma])$ as:

\begin{equation}
\vc[f]_i (\vc[\gamma]) = \begin{cases} \frac{1}{2} (\gamma_{i,1}^2 + \gamma_{i,2}^2 - \mu_i^2 \gamma_{i,n}^2) & i = 1,\dots,N_c \\
-\gamma_{i-N_c,n} & i=N_c+1,\dots,2 N_c \end{cases}
\end{equation} 

the Newton equations that define the PDIP step thus require the solution of a linear system of the form: 

\begin{equation} 
\bbm \vc[N] + \vc[\hat{M}] & \vc[B] \\ \vc[C] & \vc[E] \ebm \bbm \Delta \gamma \\ \Delta \lambda \ebm = 
\bbm r_\gamma \\ r_\lambda \ebm \label{eq:qp_Newton_system}
\end{equation} 

where $\vc[\hat{M}]$ is a diagonal matrix defined by $\vc[\hat{M}] = \sum_{i=1}^{2 N_c} \vc[\lambda]_i \nabla^2 \vc[f]_i (\vc[\gamma] ) = diag(\vc[\hat{m}])$, with $\vc[\hat{m}] =$  $[ \mu_1^2 \lambda_1,\lambda_1,\lambda_1,$ $\dots, \mu_1^2 \lambda_{N_c},\lambda_{N_c},\lambda_{N_c} ]$. $\vc[B], \vc[C]$ are banded rectangular matrices, and $\vc[E]$ is diagonal. A typical sparsity pattern for a multibody dynamics problem is shown in Fig. \ref{fig:SpyNewtonMat}. In order to find the Newton step, we must then solve this sparse linear system. We may either proceed directly, or by eliminating $\Delta \lambda$, a reduced, symmetric positive definite Schur complement matrix of size $3 N_c \times 3 N_c$ for $\Delta \gamma$ can be obtained. 

\section{The Tensor Train solver}
\label{sec:solver}
In the context of multibody dynamics, we know the system matrices for the Newton step in \eqref{eq:qp_Newton_system} to be sparse and highly structured, and dependent on the iterates $(\vc[\gamma],\vc[\lambda])$ as the PDIP iteration proceeds to the solution of problem \eqref{eq:qp}. Further, from one time step to the next, we expect the matrix required to obtain the Newton step to change in size as the active set of constraints (corresponding to pairs of objects in contact) evolves. It is because of these changes both within and between timesteps that producing an efficient and robust solution technique for the Newton system remains challenging. 

Amongst currently available hierarchical compression techniques, the \TT\ decomposition features compression and inversion algorithms that are applicable to a large set of structured matrices and that lend themselves to inexpensive global updates. In addition, they have shown to achieve sublinear precomputation times. We thus propose to use it as a framework for direct solution and preconditioning of iterative solvers for the linear systems in each PDIP iteration.  

In this section, we give a cursory description of the quantized tensor train (QTT) decomposition as a method to efficiently compress, invert and perform fast arithmetic with approximate representations of structured matrices. We then present a general discussion of its application to solving the linear systems associated with the PDIP for the CCP. This constitutes, to our knowledge, the first application of hierarchical compression solvers to the acceleration of second order optimization methods. Although it is beyond the scope of this work, we expect the techniques laid out in this section to be readily applicable to a more general class of interior point and other Newton and Quasi Newton based methods for smooth convex problems.   

\subsection{The tensor train decomposition \label{ssc:TT}}

The \TTrain\ (\TT) decomposition provides a powerful tensor compression technique \cite{oseledets2010tt} via low rank representation akin to that of a generalized SVD. We will focus on its application in the approximation of \emph{tensorized} vectors and matrices given a hierarchical subdivision of their indices, known as quantized tensor train (QTT). Matrices in this setting are further interpreted as tensorized operators acting on such tensorized vectors. We outline how this interpretation allows us to effectively employ the \TT\ as a tool for hierarchical compression and inversion of structured matrices. 

Since this approach involves considering reshaping, tensorization and vectorization of arrays, we introduce the following index notation: For the tensor multi-index $(i_1,\ldots,i_d)$, the one-dimensional index obtained by lexicographic ordering of multi-indices will be denoted by placing a bar on top, removing commas between indices: $i =\overline{i_1i_2 \cdots i_d}$. This mapping from multi-indices to one-dimensional index corresponds to the conversion of a multidimensional array to a vector which we denote $\la{b} =\vec(\cA[b])$, with $\la{b}(\overline{i_1i_2 \cdots i_d}) = \cA[b](i_1,i_2, \ldots, i_d)$. \\

\paragraph{\textbf{A motivating example}} Suppose we wish to compress the 3-tensor $\tensor{A}(i_1,i_2,i_3)$ obtained from sampling the function $f(x,y,z) = sin(x+y+z)$ on a uniform grid with $n^3$ points $(x_{i_1},y_{i_2},z_{i_3})$ in $[0,1]^3$. We can use addition formulas to decompose $f(x,y,z)$ as a sum of separable functions: 

\begin{align} 
sin(x+y+z) &= \bbm sin(x) & cos(x) \ebm \bbm cos(y+z) \\ sin(y+z) \ebm \label{eq:sine_add1} \\
                 &= \bbm sin(x) & cos(x) \ebm \bbm cos(y) & -sin(y) \\ sin(y) & cos(y) \ebm \bbm cos(z) \\ sin(z) \ebm \label{eq:sine_add2} 
\end{align} 

\begin{figure}[!t]
  \centering
  \includegraphics[width=0.7\textwidth]{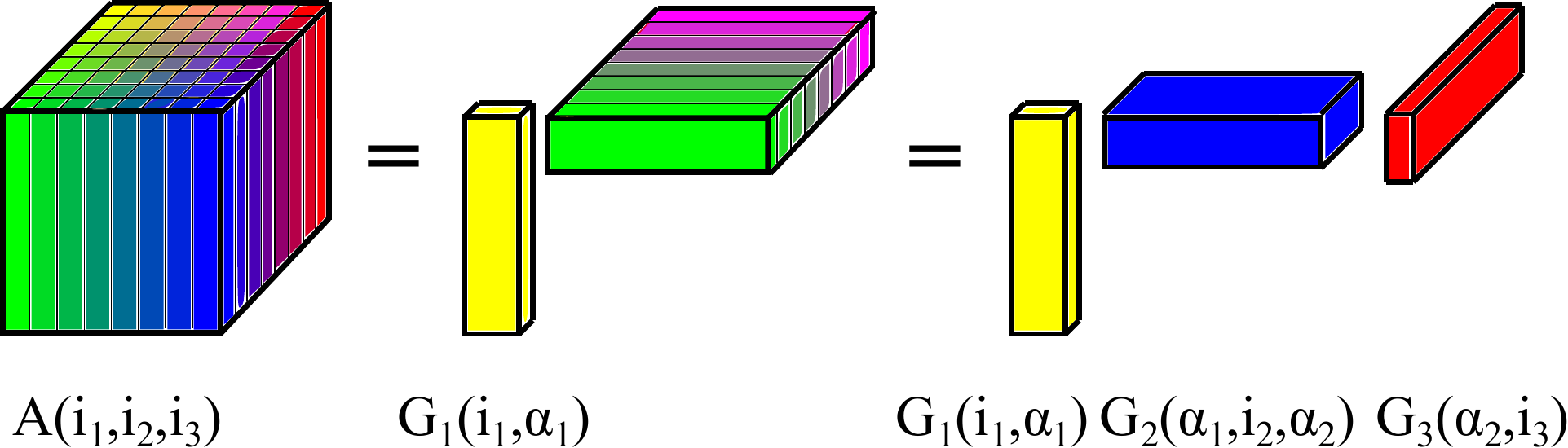}
  \caption{ \textbf{TT decomposition example: } compression of $3$ tensor from samples of $f(x,y,z)=sin(x+y+z)$ through a chain of two rank $2$ matrix decompositions corresponding to angle addition formulas in \eqref{eq:sine_add1} and \eqref{eq:sine_add2}. TT cores $\tensor{G}_i$ are obtained through a chain of low rank matrix decompositions of unfolding matrices. }
  \label{fig:TT_sinecube}
\end{figure}

Evaluation of \eqref{eq:sine_add1} and \eqref{eq:sine_add2} on the tensor $\tensor{A}$ is depicted in Fig \ref{fig:TT_sinecube}. The first step \eqref{eq:sine_add1} produces a rank $2$ decomposition $\tensor{A}(i_1,i_2,i_3) = \tensor{G}_1 (i_1,\alpha_1) \tensor{V}_1(\alpha_1,i_2,i_3)$. In \eqref{eq:sine_add2} we then decompose $\tensor{V}_1$ to separate dependency of $y$ and $z$; each of the two terms results in a rank $2$ decomposition $\tensor{V}_1(\alpha_1,i_2,i_3) = \tensor{G}_2 (\alpha_1,i_2,\alpha_2) \tensor{G}_3 (\alpha_3,i_3)$. We note that each term $G_k$ (depicted in solid yellow, blue and red in Fig \ref{fig:TT_sinecube}) depends on one dimension $i_k$ of the original tensor $\tensor{A}$, and storage has been reduced from $n^3$ to $2n+4n+2n = 8n$. 

What we have presented in this example is an exact tensor train decomposition of $\tensor{A}$. We present the analogous definition for the TT decomposition of a $d$ dimensional tensor. 

\begin{definition}
For a $d$-dimensional tensor $\cA[A](i_1,i_2,\ldots,i_d), i_k\le n_k$, sampled at $N =\prod_{k=1}^d n_k$ points, a \TT\ decomposition is of the form
\begin{equation}
  \ttA[A](i_1,i_2,\ldots,i_d) \defeq
  \sum_{\alpha_1,\ldots,\alpha_{d-1}} {\tensor{G}_1(i_1,\alpha_1)
    \tensor{G}_2(\alpha_1,i_2,\alpha_2) \dots
    \tensor{G}_d(\alpha_{d-1},i_d)},
  \label{eq:TTdeccores}
\end{equation}
where, each two- or three-dimensional $\tensor{G}_k$ is known as a
tensor core.
The ranges of auxiliary indices $\alpha_k=1, \ldots, r_k$ determine the number of terms in the decomposition. We refer to $r_k$ as the \ordinal{k} \TT-rank, analogous to matrix numerical rank. 
\end{definition}

In order to understand the chain of ``low rank'' decompositions in the tensor train format in terms of matrix low rank decompositions, we introduce auxiliary objects known as unfolding matrices. 

\begin{definition}
For a tensor of dimension $d$, the \ordinal{k} \emph{unfolding matrix}
is defined as
\begin{equation}
  \laT{A}^{k}(p_k,q_k) = \laT{A}^{k}(\overline{i_1i_2\cdots i_k} ,
  \overline{i_{k+1}\cdots i_d}) = \cA[A](i_1,i_2,\cdots,i_d)
  \quad\text{for}\quad k=1,\ldots, d,
\end{equation}
where $p_k = \overline{i_1\cdots i_k}$ and $q_k =
\overline{i_{k+1}\cdots i_d}$ are two flattened indices.
Using Matlab's notation
\begin{equation}
  \laT{A}^{k} = \reshape\left(\cA[A],\prod_{\ell=1}^k
      {n_\ell},\prod_{\ell=k+1}^d {n_\ell}\right).
\end{equation}
\end{definition}

The ranks of the \TT\ decomposition are thus the ranks of unfolding matrices. For instance, if we consider the first unfolding matrix $\laT{A}^{1}$ of tensor $\tensor{A}$ in our example, its rows will depend on $x_{i_1}$, and its columns on $y_{i_2},z_{i_3}$. Eq \eqref{eq:sine_add1} clearly implies that $\laT{A}^{1}$ is exactly of column rank $2$, and the resulting matrix decomposition

\begin{equation} \laT{A}^{1} = \tensor{G}_1 (i_1,\alpha_1) \tensor{V}_1(\alpha_1,\overline{i_2 i_3})  \end{equation} 

is interchangeable with the first rank 2 decomposition of $\tensor{A}$, requiring only to merge indices $i_2,i_3$ in $\tensor{V}_1$. We may similarly show that $\laT{A}^{2}$ is of rank $2$ (by applying addition formulas for $(x+y)$ and $z$); however, the second rank $2$ approximation in the chain is obtained by decomposing an unfolding matrix of $\tensor{V}_1$: 

\begin{equation} \laT{V}_1^{2}(\overline{\alpha_1 i_2},i_3) = \tensor{G}_2 (\overline{\alpha_1 i_2},\alpha_2) \tensor{G}_3(\alpha_2, i_3)  \end{equation}

As is the case with low rank matrix decompositions, most often a tensor $\tensor{A}$ of interest will be approximately of low rank, in the sense that given a target accuracy $\acc$, a TT decomposition with low TT ranks $\ttA[A]$ may be found such that $||\tensor{A} - \ttA[A]||_F < \acc$. This decomposition can be obtained by a sequence of low-rank approximations to $\laT{A}^{k}$. A generic algorithm proceeds as in Algorithm \ref{alg:TT-dec}.

\begin{algorithm}[!h]
  \centering
  \begin{algorithmic}[1]
    \setlength\commLen{.5\linewidth} 
    \REQUIRE Tensor $\tensor{A}$, and target accuracy $\acc$
    \STATE $\laT{M}_1 = \laT{A}^{1}$ \COMMENT{First unfolding matrix}
    \STATE $r_0 = 1$
    \FOR{$k=1$ to $d-1$}
    \STATE $[\laT{U}_k$,$\laT{V}_k] = \mathtt{lowrank\_approximation}(\laT{M}_k,\acc / \sqrt{d-1})$
    \STATE $r_k =\mathtt{size}(\laT{U}_k,2)$ \COMMENT{\ordinal{k} \TT\ rank}
    \STATE $\tensor{G}_k = \reshape\big( \laT{U}_k , [r_{k-1},n_k,r_k]\big)$ \COMMENT{\ordinal{k} \TT\ core}
    \STATE $\laT{M}_{k+1} = \reshape\Big( \laT{V}_k , [r_kn_{k+1},\prod_{\ell=k+2}^{d} {n_{\ell}}]\Big)$ \COMMENT{$\laT{M}_{k+1}$ corresponds to the \ordinal{(k+1)} unfolding matrix of $\tensor{A}$}
    \ENDFOR
    \STATE $\tensor{G}_d = \reshape\big( \laT{M}_d, [r_{d-1},n_d,1]\big)$ \COMMENT{Set last core to the right factor in the low rank decomposition}
    \RETURN $\ttA[A]$
  \end{algorithmic}
  \caption{\TT\ decomposition}
  \label{alg:TT-dec}
\end{algorithm}

Due to the cost of successive low rank factorizations, implementing Algorithm \ref{alg:TT-dec} will predictably lead to relatively high computational cost, $O(N)$ or higher, which is exponential in the tensor dimension $d$. Instead, we employ TT rank revealing strategies based on the multi-pass AMEN (Alternating Minimal Energy) Cross algorithm of \cite{oseledets2010tt}, in which a low-\TT-rank approximation is initially computed with fixed ranks and is improved upon by a series of passes through all cores. 

\paragraph{\textbf{\TT\ compression update}} We note that, since these algorithms obtain the \TT\ approximation based on an iterative, greedy rank detection procedure, they allow for inexpensive updates to a \TT\ compressed representation of a tensor $\tensor{A}$. If we wish to produce the \TT\ representation of $\tensor{\tilde{A}} = \tensor{A} + \tensor{E}$, the AMENCross algorithm may be started using the \TT\ representation of $\tensor{A}$ as an initial guess. If $\tensor{E}$ has small \TT\ ranks, this provides a significant speed-up for this compression algorithm, which often converges in a few iterations to the updated representation. 

\paragraph{\textbf{Computational complexity and memory requirements}} Because AMEN Cross and related TT rank revealing approaches proceed by enriching low-\TT-rank approximations, all computations are performed on matrices of size $r_{k-1}n_k \times r_k$ or less. In \cite{corona2015tensor}, it is shown that  the complexity for this algorithm is thus bounded by $O(r^3 d)$ or equivalently $O(r^3 \log N)$, where $r = \max(r_k)$ is the maximal \TT-rank that may be a function of sample size $N$ and accuracy $\acc$. For a large number of structured matrices as well as their inverses, $r$ typically stays constant or grows logarithmically with $N$ \cite{kazeev2012low, oseledets2010approximation, oseledets2011tensor, corona2015tensor}.The overall complexity of computations is then \emph{sublinear} in $N$.

\begin{figure}[!hb]
  \centering
  \includegraphics[width=\textwidth]{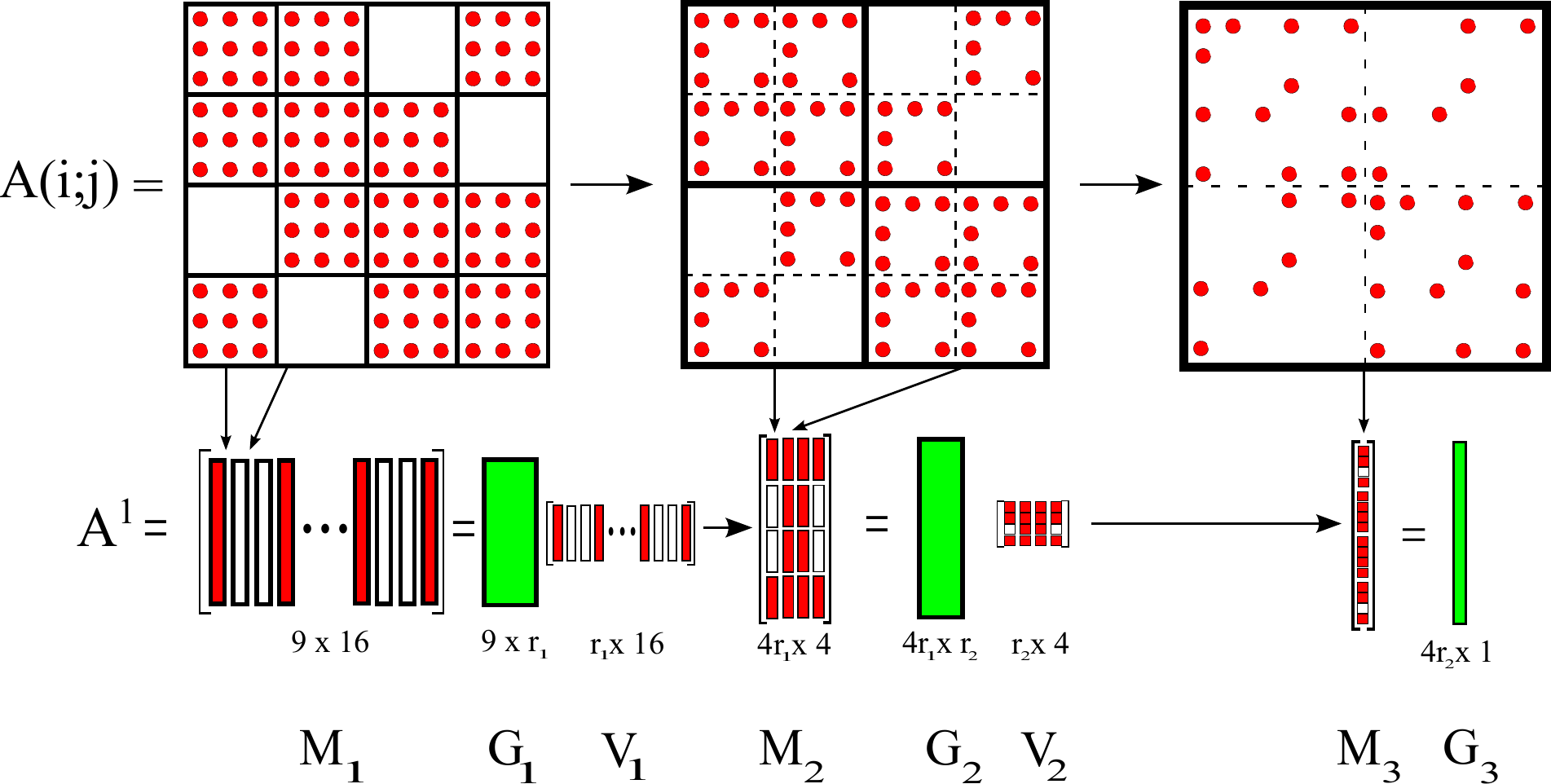}
  \caption{ \textbf{TT decomposition for a matrix with $d=3$.} Each level of refinement for the matrix block hierarchy corresponds to one tensor dimension and TT core. Columns of the corresponding unfolding matrix are obtained by vectorizing matrix blocks (in red if non-zero, white otherwise). The example illustrated above uses an interpolatory low rank decomposition, producing a uniform subsampling of tree nodes. The TT decomposition process thus consists of finding a hierarchical basis of matrix block entries. We show the main steps in algorithm \ref{alg:TT-dec}: (i) computing a low rank decomposition $\laT{U}_k \laT{V}_k$  for the corresponding unfolding matrix $\laT{M}_k$; (ii) taking $\laT{U}_k$ (in green) as the \ordinal{k} TT core $\tensor{G}_k$;  and (iii) interpreting the right factor $\la{V}_k$ as a matrix in level $k+1$ to form $\laT{M}_{k+1}$. 
  }
  \label{fig:TT_comp}
\end{figure}

\subsection{\TT\ for hierarchically structured matrices}\label{ssc:TTmat}
Consider a block-sparse, structured matrix $\laT{A}$. For simplicity of presentation we assume matrix rows and columns to be of size $N = n_1 2^{d-1}$. We then recursively bisect them, representing the resulting hierarchy with a binary tree $\tree{T}$ (number of children $n_k=2$) with $n_1$ points in each leaf node. Let $\tree{T}_\lbl{\Pi}$ then denote the product tree $\tree{T} \times \tree{T}$; nodes on this tree correspond to pairs of source and target nodes. $\tree{T}_{\Pi}$ can be understood as a hierarchy of matrix-blocks, or equivalently, of all interaction operators between subsets $B_i$ and $B_j$ at a given level $\ell$ of $\tree{T}$. Each node on this tree can thus be indexed by integer row and column index coordinate pairs $(i_{k},j_{k})$ with $k\le \ell$. Equivalently, we can consider block integer coordinates $b_{k} \in
\{1,\cdots,n^2_{k}\}$ for $\tree{T}$ such that $b_{k} = \overline{i_{k}j_{k}}$. We then apply the \TT\ decomposition to the corresponding tensorized form of $\laT{A}$, $\tensor{A}_{\tree{T}}$, a \mbox{$d$-dimensional} tensor with entries defined as

\begin{equation}
  \tensor{A}_{\tree{T}}(b_1,b_2,\ldots,b_d) =
  \tensor{A}_{\tree{T}}(\overline{i_1j_1},
  \overline{i_2j_2},\ldots,\overline{i_dj_d}) =
  \laT{A}(\overline{i_1i_2 \cdots i_d}, \overline{j_1j_2 \cdots j_d}),
\end{equation}

and obtain an approximate \TT\ factorization $\ttA[A]$. Each core of $\ttA[A]$, $\tensor{G}_{k}(\alpha_{k-1}, \overline{i_kj_k}, \alpha_{k})$ depends only on the pair of source and target tree indices at the corresponding level of the hierarchy. When performing matrix arithmetic, such as matrix-vector product or inversion, TT cores are often reshaped as
$n_k \times n_k$ matrices parametrized by $\alpha_{k-1}$ and $\alpha_k$.

In Fig. \ref{fig:TT_comp}, we demonstrate the \TT\ decomposition algorithm applied to a matrix $\laT{A}$ for binary source and target trees with depth $d=3$ and three points in the leaf nodes $n_1 = m_1 = 3$, implying $N=12$. In this example, the tree $\tree{T}_\lbl{\Pi}$ is a matrix-block quadtree. At each level of the hierarchy, every column of the corresponding unfolding matrix $\la{A}^k$ is a vectorized form of the block $A_{ij}$. We note that, for a block-sparse matrix $\laT{A}$, the majority of these blocks, and hence columns of $\la{A}^k$ at intermediate levels $k$ of the tree are identically zero. This, along with any additional block structure in $A$ ensures that $A^k$ is approximately low rank.   

\paragraph{\textbf{Fast arithmetic and TT direct solvers}}

Fast methods for TT matrix arithmetic are available for a number of operations, including inversion and matrix-vector and matrix-matrix products. The TT solvers in this work find an approximate TT structure for the inverse and employ the corresponding TT matrix-vector product $\gamma = \la{A}^{-1} b$. This mat-vec algorithm proceeds by contracting one dimension of the tensorized matrix $A^{-1}$ at a time, applying the \ordinal{k} \TT\ core. Its complexity can be shown to be $O(r^2 N \log N)$. 

\TT\ inversion methods compute an approximate \TT\ decomposition of $\laT{A}^{-1}$ given the \TT\ decomposition of $\laT{A}$. In these algorithms, matrix equations for each tensor core of the inverse is solved iteratively. Following an Alternating Least Squares (ALS) algorithm, given an initial guess for the inverse in \TT\ form, it proceeds by iteratively cycling through the cores (freezing all cores but one) and solving a linear system to update the \ordinal{k} core of $\laT{A}^{-1}$. \TT\ ranks of the inverse are not known a priori (and are distinct to those for $\laT{A}$), and so strategies to increase core ranks are needed to ensure convergence of the ALS procedure to an accurate inverse representation. Further details about variants of this approach can be found in \cite{oseledets2012solution, dolgov2013amen1, dolgov2013amen2}. The complexity of this algorithm for a maximum TT rank $r$ for both $\ttA[A]$ and $\ttA[A]^{-1}$ is bound by $O(r^4 \log N)$, as shown in \cite{corona2015tensor}. 

\paragraph{\textbf{Inverse compression update}} We note that, since they share the same iterative, greedy rank detection structure with the AMENCross compression algorithm, \TT\ direct solver routines may also be significantly sped up using an approximate initial solution. In the context of the Newton step in interior point methods, we expect both forward and inverse operators to be obtainable via low \TT\ rank updates.   

\subsection{TT Newton system solver}\label{ssc:TTNewton}

We now discuss how to adapt the tensor train decomposition framework to accelerate computation of the Newton steps within the PDIP method applied to the CCP. We first show how at a given timestep, the tensor train provides an easy-to-update approximate inverse for the Newton system's Schur complement matrix. We then describe a procedure to hot-start the TT compression and inversion algorithms re-using information from the previous timestep, even when the corresponding set of contacts (and thus, matrix size and structure) change. 

\paragraph{\textbf{PDIP iteration}} At a given timestep $t$, the \ordinal{k} PDIP iteration involves the solution of the Newton system in \ref{eq:qp_Newton_system}, with a $5N_c \times 5N_c$ system matrix of the form

\begin{equation} \vc[A]_k = \bbm \vc[N] + \vc[\hat{M}]_k & \vc[B]_k \\ \vc[C]_k & \vc[E]_k \ebm \end{equation}

where $\vc[N]$ is fixed and the rest of the matrix blocks involved depend on the current iterate $\vc[\gamma]_k, \vc[\lambda]_k$. We may, alternatively, choose to solve the Schur-complement system for $\Delta \vc[\gamma]$, resulting in the $3N_c \times 3 N_c$ matrix: 

\begin{equation} \vc[S]_k = \vc[N] + \vc[\hat{M}]_k -\vc[B]_k\vc[E]_k^{-1}\vc[C]_k \end{equation} 

We note that working with the Schur complement $\vc[S]_k$ is generally preferable, as it is symmetric positive definite and smaller in size. We observe that its sparsity pattern across iterations is always that of matrix $\vc[N]$. That is due to the fact that $\vc[\hat{M}]_k$ is diagonal and $\vc[B]_k\vc[E]_k^{-1}\vc[C]_k$ is block-diagonal ($3 \times 3$ blocks).  

\begin{figure}[!h]
  \centering
  \includegraphics[width=0.8\textwidth]{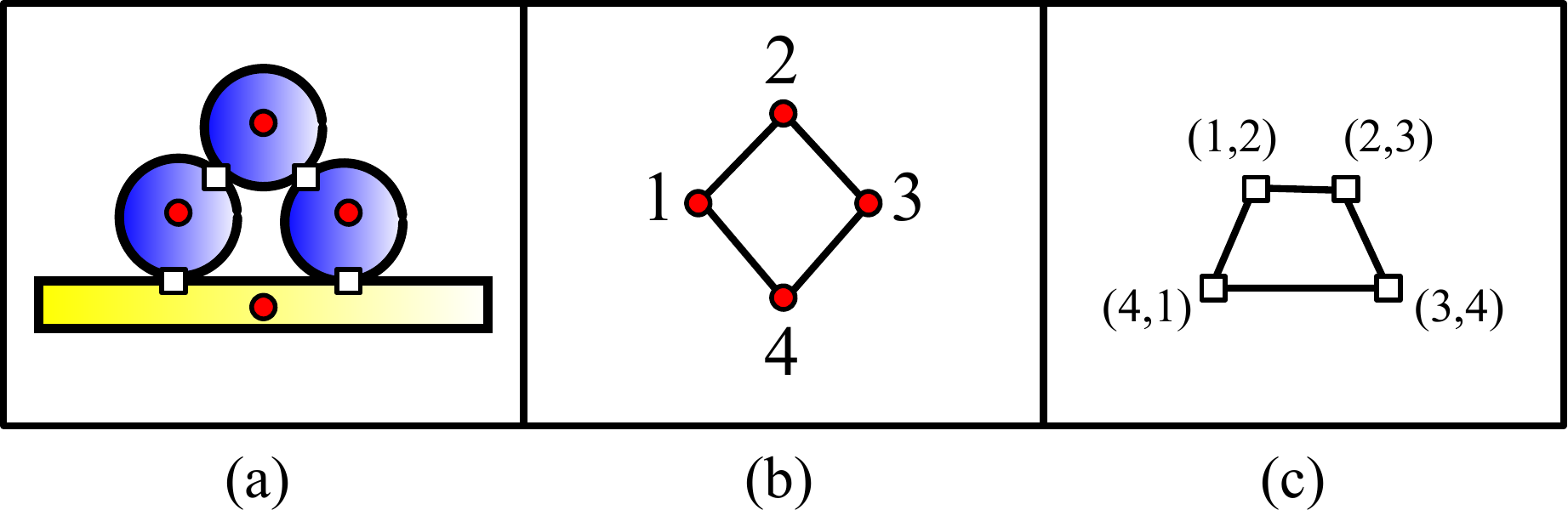}
  \caption{ \textbf{Collision graphs example: } (a) configuration of three spherical bodies lying on a flat surface. Centers of mass are depicted as red circles, collision points as white squares. (b) Graph connecting bodies that are in contact (c) Graph connecting contacts sharing a body in common. This graph reveals the block-sparsity pattern for $\vc[N]$ and $\vc[S]$. }
  \label{fig:collisiongraphs}
\end{figure}

We recall that matrix $\vc[N]$, the Hessian of the quadratic $q(\vc[\gamma])$, is of the form $\vc[D]^T \vc[M]^{-1} \vc[D]$ (Eq \eqref{eq:Ndef}), with $\vc[D] \in \mathbb{R}^{6M \times 3N_c}$ the contact transformation matrix and $\vc[M]$ the diagonal mass matrix. In order to understand the sparsity pattern of $\vc[N]$, we partition it into $3 \times 3$ blocks corresponding to each contact. The $(\ell,i)$th block of $\vc[D]$ is non-zero if the \ordinal{i} contact involves body $B_\ell$. As a result, the $(i,j)$th block of $\vc[N]$ is non-zero if the \ordinal{i} and \ordinal{j} contacts share a body in common. In Fig. \ref{fig:collisiongraphs}, we show a simple example with three spherical bodies lying on a flat surface. Given the network of bodies at contact in (b), the sparsity pattern in $\vc[N]$ corresponds to the edge-adjacency graph in (c); two two edge nodes (indexed by body pairs) are connected if they share a body in common. 

A hierarchy for contact pairs may be generally constructed from the corresponding adjacency graph using graph partition techniques such as the nested dissection method. For the cases of interest in this work, this hierarchy can be readily obtained from a spatial octree hierarchy of the positions of contact pairs in 3$d$ space, as long as care is put to separate those associated with large contact geometries (e.g. those associated with walls or containers such as the flat surface in Fig \ref{fig:collisiongraphs}).

As indicated in Section \ref{ssc:TTmat}, given a hierarchical partition of matrix indices, at each timestep we tensorize the Schur complement $\vc[S]_k$, and for a given target accuracy $\acc$ we construct approximate TT decompositions $\ttA[S]_k$ and $\ttA[S]_k^{-1}$. We then have the option to use $\ttA[S]_k^{-1}$ as a direct solver, or couple it with an iterative procedure if a more accurate solve is needed (e.g. iterative refinement, or a preconditioned Krylov subspace method). We then consider $S_{k+1}$ as a perturbation 

\begin{equation}
S_{k+1} = S_k + L_k,
\end{equation}

where $L_k$ is block-diagonal. While $L_k$ is generally of matrix rank $O(N_c)$, across all experiments in Section \ref{sec:results} we observe it to be of approximate low TT rank, and the same is found for $S_{k+1}^{-1}$ as a perturbation of $S_{k}^{-1}$. This fact allows us to use the TT decompositions $\ttA[S]_k,\ttA[S]_k^{-1}$ to hot-start the corresponding compression and inversion algorithms, reducing precomputation times considerably. 

\paragraph{\textbf{Re-using information across timesteps}} Employing information from the solution of the CCP at a timestep $t$ to hot-start the PDIP iteration at the next timestep $t+\Delta t$ is notoriously hard; even if it can be used to produce a feasible point that is close to the optimum (which is non-trivial due to changes in the set of contacts), efforts by the PDIP algorithm to preserve centrality might cause it to take small steps and waste time ``returning'' to the central path.

For this reason, we initialize each PDIP iteration by making the tangential force impulses $\gamma_{i,1},\gamma_{i,2}$ equal to zero, and the normal force $\gamma_{i,n}$ equal to either a constant preset value (e.g. $1$) or to the value computed in the previous timestep if the corresponding contact persists across timesteps. This ensures that our initial value is feasible and lies safely inside the friction cone. 

\begin{figure}[!h]
  \centering
  \includegraphics[width=0.5\textwidth]{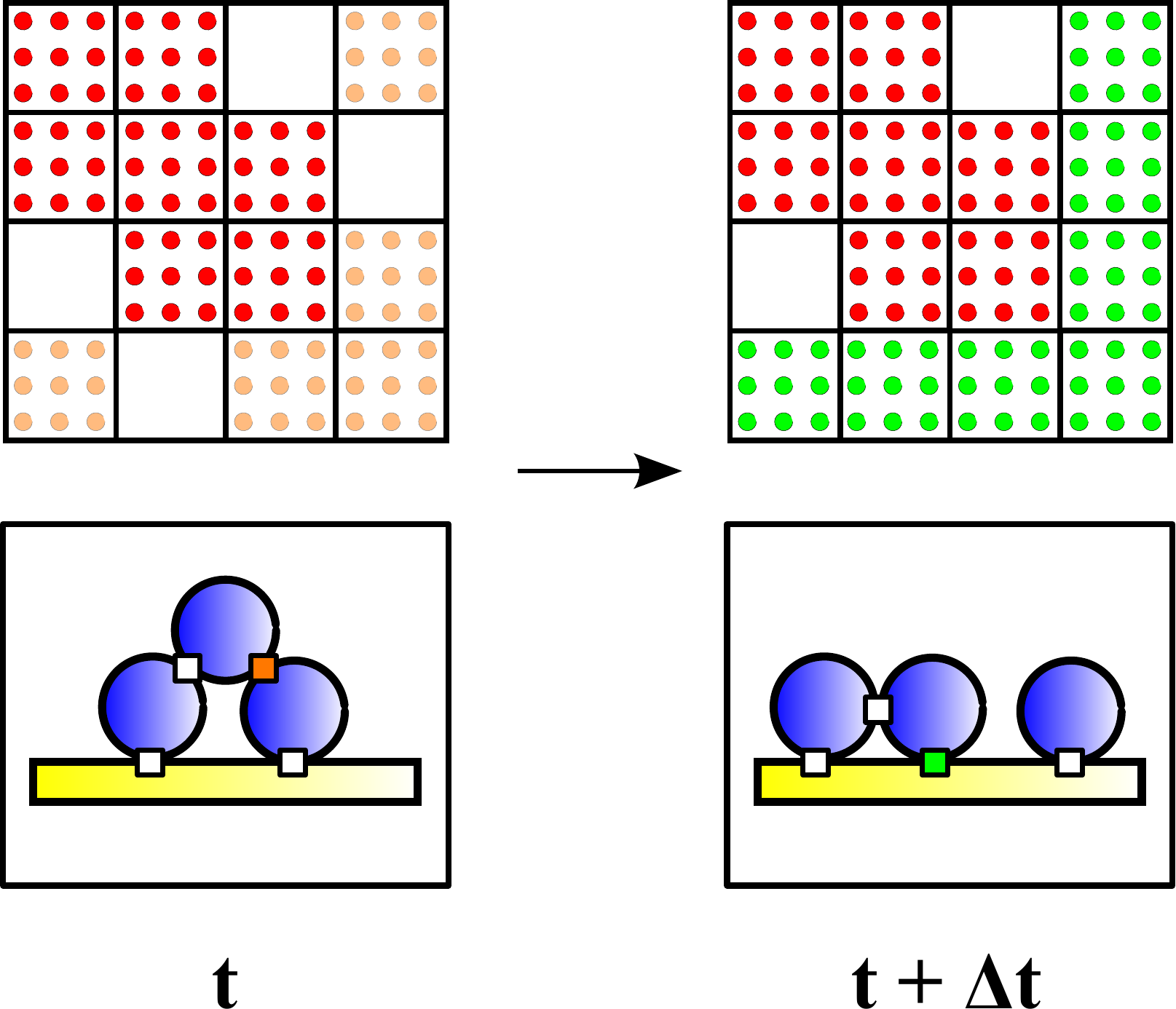}
  \caption{ \textbf{Change in contact set and matrix entries: } From time $t$ to $t + \Delta t$, we depict a change in the configuration of the three bodies in Fig \ref{fig:collisiongraphs} as the system evolves. Above each configuration, we depict the corresponding block-sparsity pattern for the initial Schur complement matrices. Persistent matrix entries are depicted in red, entries being removed in orange and entries being introduced in green. }
  \label{fig:changecontactset}
\end{figure}

Since pairs of bodies may phase in and out of contact, the set of contacts considered at each timestep changes. However, as long as $\Delta t$ and relative velocities are sufficiently small, it is likely that the set of persistent contacts from one timestep to the next will be large. In all validation experiments in Section \ref{sec:results}, in fact, well over $90\%$ of contacts persist once objects have sedimented. We may then re-use the TT decompositions for the \emph{initial} Schur complement matrix. Care must be taken to make this decomposition compatible, introducing new contacts into the hierarchy and ``deleting'' contacts that have ceased to exist. A simple example of this is depicted in Fig \ref{fig:changecontactset}.   

In order to use the TT-factorization at timestep $t$ as shown in Fig \ref{fig:TT_comp} to form an initial approximation for compression and inversion at the next timestep, we need to reconfigure the TT cores so that they remain an approximation for the submatrix for persistent contacts (in red). For instance, if they correspond to a hierarchical interpolation of matrix entries, we could keep the first $d-1$ cores (containing interpolation weights) and update the entries in $\tensor{G}_d$ corresponding to new contacts (in green).  

\begin{figure}[!h]
  \centering
  \includegraphics[width=0.5\textwidth]{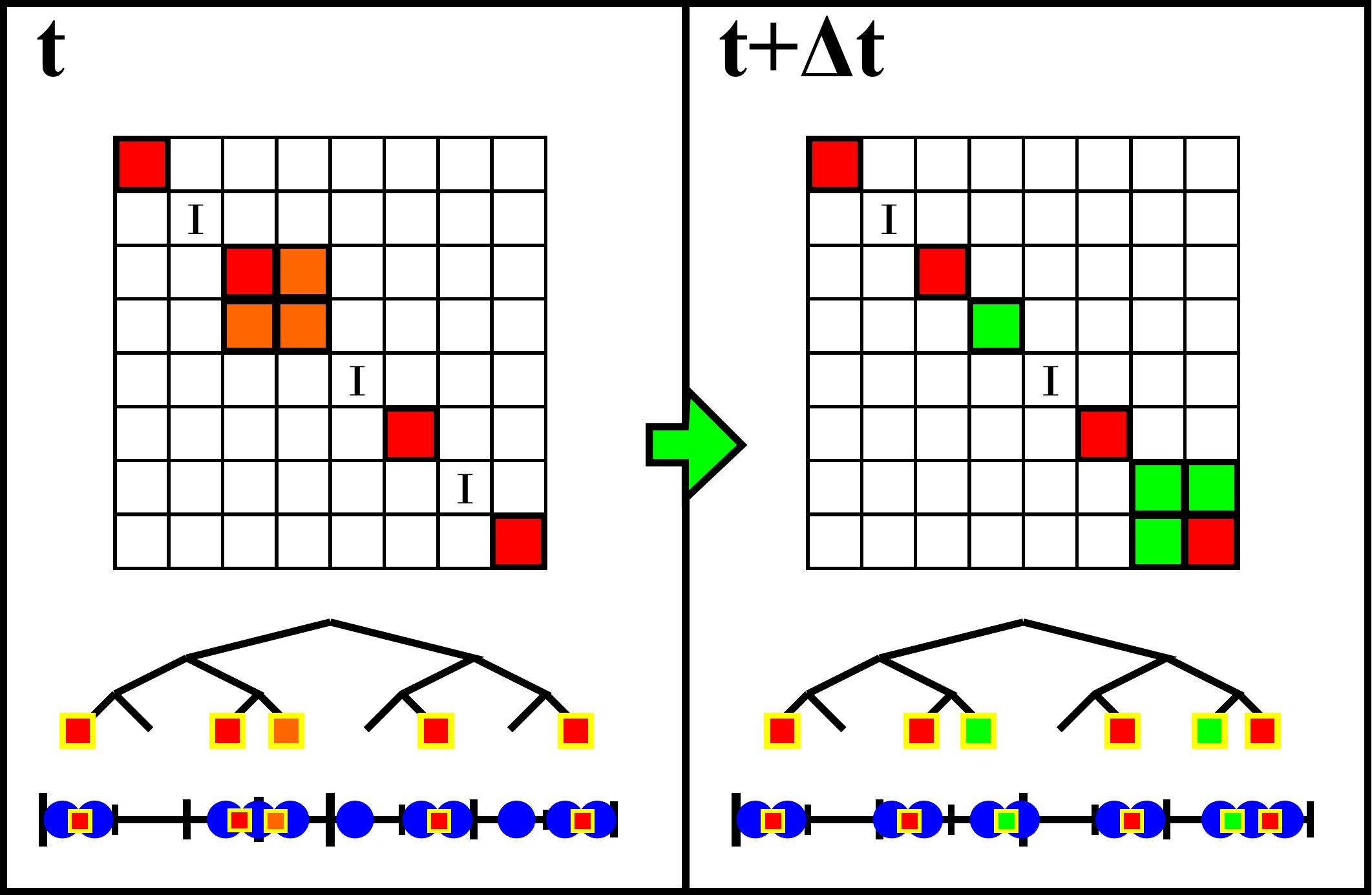}
  \caption{ \textbf{TT tensor index hierarchy evolution: } We depict the evolution of a spatial hierarchy used to index pairs of bodies at collision from a simple one-dimensional example as it evolves from timestep $t$ to $t + \Delta t$. Above each configuration, we show the corresponding binary tree and sparse matrix structure for the Schur complement matrix. Persistent collision pairs and corresponding matrix entries are depicted in red, pairs being removed in orange and pairs being introduced in green. In this example, one of 5 nodes is removed and two new nodes are introduced. }
  \label{fig:hierarchyevolve}
\end{figure}

We note that assigning tensor indices to ``new body pairs'' and eliminating ``old'' ones requires us to track and modify the spatial hierarchy that is used to encode our matrix as a tensorized array. This requires a certain degree of flexibility and adaptivity, if representations for $t$ and $t+\Delta t$ are to be compatible. In this work, we resolve this issue by considering our hierarchy as an adaptive tree which is itself a subset of a uniform tree with $\tilde{N}_c = 2^L$ elements. Any elements outside our current hierarchy are dummy variables, and the corresponding matrix for this augmented set of degrees of freedom is a permutation of a matrix with two diagonal blocks, $A$ and $I$ the identity. We detail a small one-dimensional example of this in Fig \ref{fig:hierarchyevolve}. 

This setup allows us to produce a compatible initial guess for $\ttA[A]$ and $\ttA[A]^{-1}$ across timesteps as long as the adaptive tree can be properly updated and still fits within the regular tree with $\tilde{N}_c$ leaf nodes. New nodes take the place of dummy variables, and persistent nodes are potentially re-indexed if they move too far. Otherwise, we reset this structure and compute the TT factorizations from scratch.   

\section{Numerical Results}
\label{sec:results}
We now demonstrate the performance of the TT-based solver when accelerating the solution of the Newton step system, and thus of second order methods such as PDIP, in the context of dense, multiple rigid-body dynamics. As mentioned in Sections \ref{sec:intro} and \ref{sec:solver}, we know that given a desired target accuracy $\varepsilon$, if the associated TT ranks $r_k$ are bounded or slowly growing as a function of problem size $N$, factorization costs and storage for the TT solver are \emph{sublinear} in $N$, and that applying this matrix to a given right-hand-side is $O(N \log N)$. We wish to study if Newton system matrices are TT compressible in this sense, and to test their performance and scaling in ther solution as problem size grows (determined by the number of collisions $N_c$). By harnessing the ability of the TT to produce economic updates both from one Newton iteration to the next, as well as across timesteps, we demonstrate significant improvement for the solution of the complementarity problem using second order methods. 

The complementarity method for frictional contact, as well as the solution methods for its CCP relaxation have been thoroughly validated and contrasted with experimental data \cite{melanz2016experimental,melanz2017comparison}. We set up three experiments based on standard phenomena in terramechanics, focusing on how TT-based linear solvers perform in the PDIP iteration. Models and contact dynamics simulations are performed using open-source library Project Chrono \cite{projectchrono}, PDIP solver comparisons are performed with a serial Matlab implementation and all TT methods are based on the Matlab TT-Toolbox \cite{toolboxtt}. All experiments are run on the serial queue of University of Michigan's Flux computing cluster. 

\subsection{TT compressibility and information re-use}

We first carry out an assessment to determine how the general behavior of TT ranks for the Schur matrix and its inverse vary with target accuracy $\varepsilon$ and maximum allowed TT ranks, varying these parameter from $10^{-2}$ to $10^{-4}$ and $r \leq 10,100,1000$, respectively. This gives us a general idea of the compressibility of these matrices upon re-ordering their nodes according to a spatial hierarchy. Since algorithmic constants for TT compression, inversion and matrix-vector apply all depend on rank, this also informs out choice of $\varepsilon$. We present average results for $1000$ PDIP iterations for a sedimentation experiment with $M = 35939$ rigid bodies in Section \ref{ssc:PDIPcompare}; we note these are largely replicated across experiments, and that ranks and compression times grow very slowly with the number of degrees of freedom. 

\begin{table}[!htb]
\centering

\begin{tabular}{ l | c  c  c |}
              & $r\leq10$ & $r\leq100$ & $r\leq1000$ \\[1 mm]
\hline
$\varepsilon=10^{-2}$ & $4.6 / 4.4$ & $10.5 / 6.7$ & $10.5 / 6.7$ \\ 
$\varepsilon=10^{-3}$ & $8.7 / 8.5$ & $59.3 / 45.2$ & $65.2 / 46.4$ \\
$\varepsilon=10^{-4}$ & $9.1 / 9.0$ & $68.6 / 58.1$ & $113.3 / 91.2$ \\
\hline
\end{tabular}
\caption{\textbf{Average TT ranks of Schur matrix:} Setting target accuracy $\varepsilon$ and maximum rank $r$, we record average ranks for $1000$ timesteps of a sedimentation simulation of $M = 35939$ rigid bodies.}
\label{tbl:avgttranks}
\end{table}

\begin{table}[!htb]
\centering

\begin{tabular}{ l | c  c  c |}
              & $r\leq10$ & $r\leq100$ & $r\leq1000$ \\[1 mm]
\hline
$\varepsilon=10^{-2}$ & $10$ & $25$ & $25$ \\ 
$\varepsilon=10^{-3}$ & $16$ & $1244$ & $1317$ \\
$\varepsilon=10^{-4}$ & $16$ & $2740$ & $8157$ \\
\hline
\end{tabular}
\caption{\textbf{Average TT compression times for Schur matrix:} Setting target accuracy $\varepsilon$ and maximum rank $r$, we record average compression times (in seconds) for the entire PDIP solver ($\sim50-100$ iterations per timestep) for $1000$ timesteps of a sedimentation simulation of $M = 35939$ rigid bodies.}
\label{tbl:avgttcomp}
\end{table}

In Table \ref{tbl:avgttranks} we observe that the Schur matrices and their inverses are indeed highly compressible, and that as we increase the maximum allowed rank, individual and average TT ranks converge. This is replicated in all our experiments, regardless of number of particles and number of collisions. In Table \ref{tbl:avgttcomp} we record the combined matrix compression and inversion times; as indicated in Section \ref{ssc:TTmat}, performance of these algorithms is bounded by terms of the form $r^k \log N_c$. Given that rank growth results in a substantial increase in precomputation times, in our experiments we find that a TT preconditioner approach with $\varepsilon = 10^{-2}$ is most effective in reducing overall computation for the PDIP iterations. We note, however, that this parameter selection is generally problem and implementation dependent. 

We then wish to quantify the impact of the strategies described in Section \ref{ssc:TTNewton} to re-use information in TT factorizations for the matrix and its inverse. For this purpose, we run two versions of the TT preconditioned PDIP iteration with and without factorization re-use for 100 timesteps for the same sedimentation problem described above. We set target accuracy to $\varepsilon = 10^{-2}$ and bound maximum ranks at $r\leq10$. 

\begin{figure}[!htb]
\centering
\includegraphics[width=\textwidth]{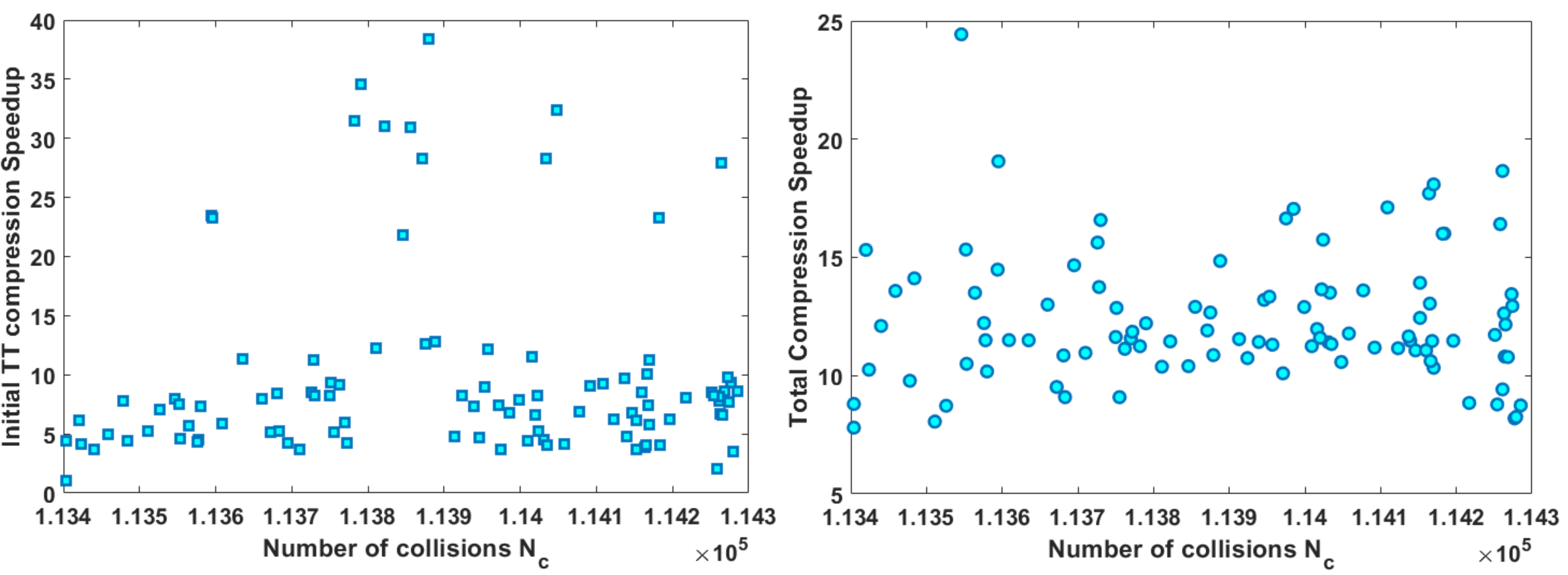}
\caption{\textbf{TT compression speedups due to information re-use:} \emph{Left:} we plot the ratio of compression times without and with information re-use for the initial Newton step for each timestep against number of collisions $N_c$; \emph{Right:} we plot the same ratio for total compression times for the PDIP iterations.}   
\label{fig:TT_comp_strategy}
\end{figure}

In Fig \ref{fig:TT_comp_strategy}, we plot precomputation times for the initial Newton steps and for the entire PDIP iteration, in order to measure how effective our information re-use strategies for the TT are in bringing down compression and inversion costs. For the initial PDIP iteration, we see that except for the first timestep (for which both methods have no prior information), information re-use always provides a speedup, between 5 and 15x for most timesteps shown here. The accumulated effect of both re-use strategies can be seen in the second plot on the right; precomputation is improved for all timesteps, resulting in a 10 to 15x overall speedup. Across all experiments, we observe both an improvement by about an order of magnitude and a drastic reduction in the variance of precomputation times, making the TT approach faster and more robust. We also note that predictably, this speedup is larger for higher target accuracies, as the iterations through all TT cores become more computationally burdensome. 

\subsection{Performance comparison experiments for PDIP iteration}\label{ssc:PDIPcompare}  

\begin{figure}[H]
\centering
\includegraphics[width=\textwidth]{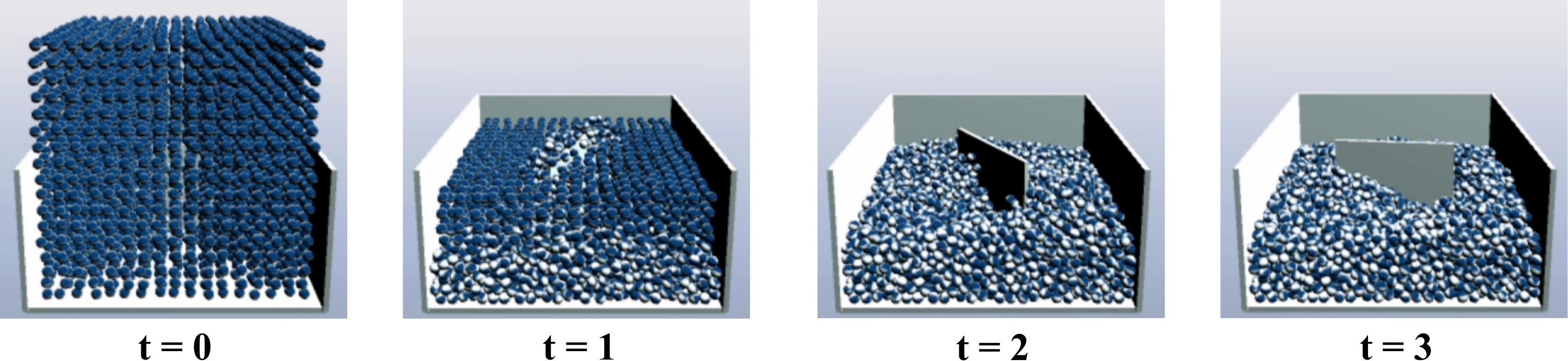}
\caption{\textbf{Sedimentation with rotating mixer:} snapshots of simulation for sedimentation of $4913$ spheres on a box shaped container with a rotating mixer with constant angular velocity. }   
\label{fig:mixer}
\end{figure}

\paragraph{Sedimentation on box with rotating mixer} A randomly pertubed cubic lattice of $(2n+1)^3$ rigid particles of spherical shape of radius $0.1$ and friction coefficient $\mu = 0.25$ are dropped and sediment under gravity into a fixed box with a rotating mixer, as shown in Fig. \ref{fig:mixer}. We then run simulations for $n=8,16,32$, with a total of rigid bodies $M = 4915,35939$ and $274627$ (including the box and mixer), scaling up the box size to keep particle density roughly constant. We set a timestep size $\Delta t = 0.025$ and target accuracy $1e$-$4$ for the PDIP solver, until such time as all objects are deposited in the box and undergoing mixing. 

We compute the Newton step through the solution of the corresponding Schur complement system. In order to test performance of the TT-based preconditioner, we compare precomputation and solution times for a preconditioned biconjugate gradient stabilized (BICGSTAB) method against Incomplete LU and unpreconditioned versions of this iterative solver. We note that while use of the conjugate gradient (CG) method may be generally preferrable for these systems, we observe its performance may degrade due to ill-conditioning as iterates approach the feasible set boundary, and so for simplicity of presentation we exclude it from our comparison.  

As discussed in Section \ref{ssc:TTNewton}, the Tensor Train approach allows us to produce approximate direct and preconditioned iterative solvers by varying target accuracy and maximum TT rank in the compression and inversion processes.   Following preliminary parametric studies, we find that setting target accuracy for TT inversion to $1e$-$2$ and capping maximum TT ranks at $r \leq 10$ provides the best performance in terms of the trade-offs involved in precomputation and TT preconditioner apply costs for our experimental setup. 

We note that for practical purposes, the maximum number of iterations for the unpreconditioned solve was set at $1000$; this limit was often reached reducing the overall accuracy of the resulting linear system solution and thus the quality of the PDIP iteration. Through further testing, we consistently observe average BICG iteration counts of $3$ - $5 \times 10^3$ and a $5$-$6$ fold increase in computational cost when removing this constraint. These estimates should be considered whenever comparing either of the preconditioned methods against unpreconditioned BICG.   

\begin{figure}[H]
\begin{center}
\includegraphics[width=0.9\textwidth]{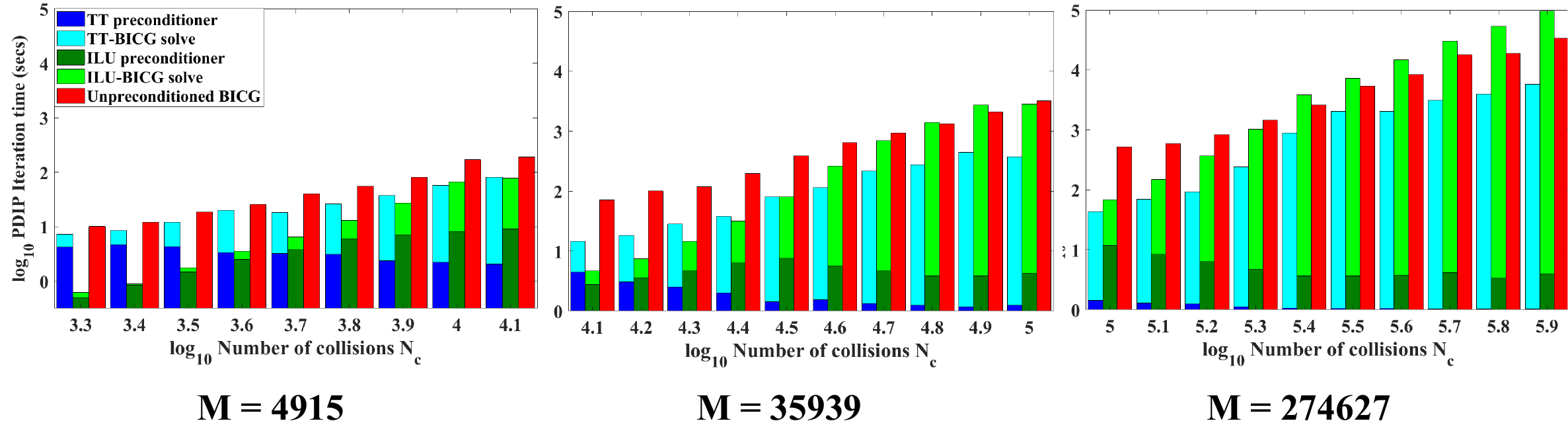}
\end{center}
\caption{\textbf{Log-log plot solver comparison:} for each experiment, we bin according to $log_{10} N_c$ and plot average PDIP iteration times; Unpreconditioned BICGSTAB in red, ILU-BICGSTAB in green and TT-BICGSTAB in blue. For preconditioned solvers, we display the proportion spent in precomputation in a darker shade, solve times in a lighter one.}
\label{fig:sedim_totstk}
\end{figure}

In Fig. \ref{fig:sedim_totstk} we can observe how performance for each of these solvers scales with the number of collisions $N_c$. Two factors are contributing to increase problem complexity in this experiment: as objects sediment in the box, the number of collisions tends to increase and the Schur complement matrix becomes less sparse. From this plot, we can readily observe that the TT preconditioned solver generally shows superior scaling, outperforming the ILU preconditioner at about $N_c \sim 20000$, and gaining an order of magnitude speed-up against the ILU preconditioner by the end of experiments with $M = 35939,274267$. We can observe that while the fraction of time spent in precomputation tends to a constant for ILU (indicating similar asymptotic scaling for precomputation and solve times, experimentally $O(N_c^3)$, it quickly goes down to zero for the TT.    

\begin{figure}[H]
\centering
\includegraphics[width=0.9\textwidth]{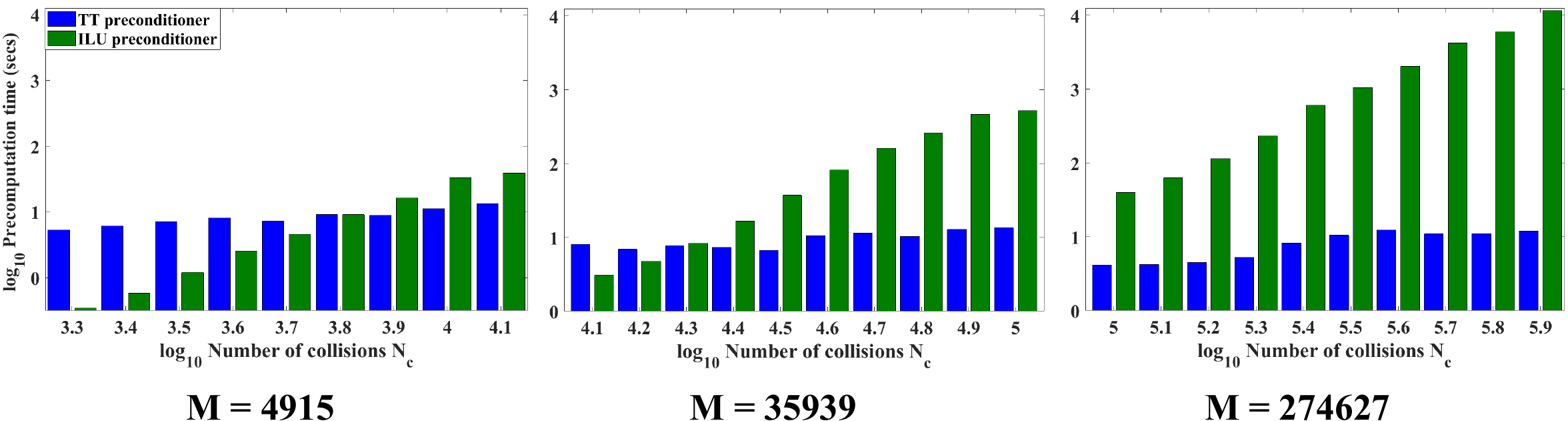}
\caption{\textbf{Log-log plot precomputation comparison:} for each experiment, we bin according to $log_{10} N_c$ and plot average preocmputation times; ILU-BICGSTAB in dark green and TT-BICGSTAB in dark blue.}
\label{fig:sedim_prec}
\end{figure}

In Fig. \ref{fig:sedim_prec}, we take a closer look at average precomputation times for each PDIP iteration for TT and ILU. We note that this involves computing roughly $50$-$100$ factorizations, one per Newton iteration. In all experiments, we confirm that compression and inversion times for the Tensor Train approach grow extremely slowly with $N_c$, staying on a range from 5 to 15 seconds. 

\begin{figure}[H]
\centering
\includegraphics[width=0.9\textwidth]{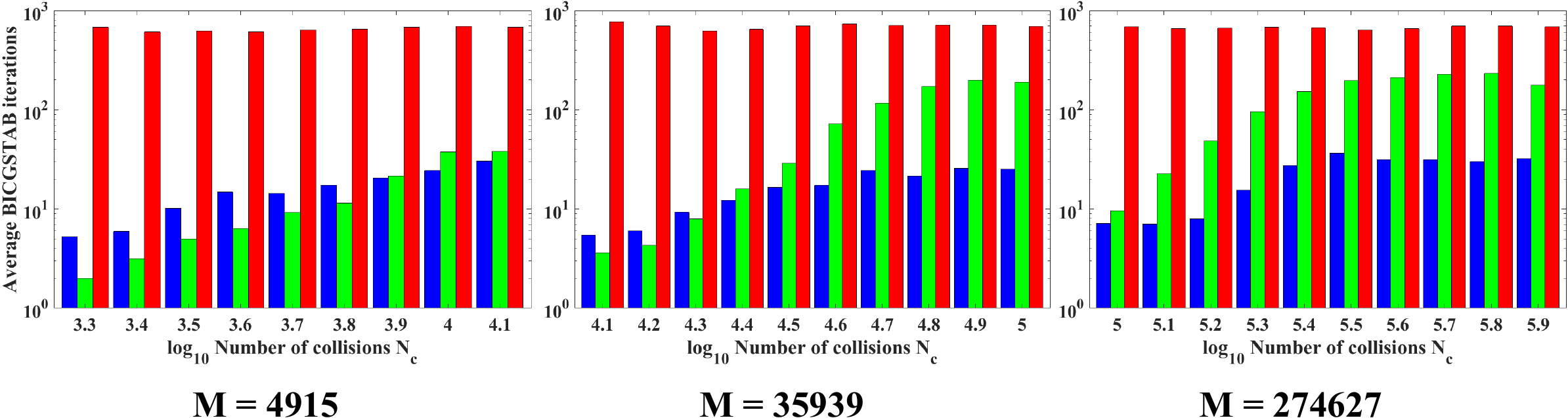}
\caption{\textbf{Average iteration count comparison:} for each experiment, we bin according to $log_{10} N_c$ and compare average BICGSTAB iteration counts; Unpreconditioned BICGSTAB in red, ILU-BICGSTAB in green and TT-BICGSTAB in blue.}
\label{fig:sedim_iter}
\end{figure}

Finally, we compare BICGSTAB average iteration counts in Fig. \ref{fig:sedim_iter}. We note that while iteration counts generally increase at the beginning of the experiment, those for the TT grow slower and settle sooner as particles sediment; for $M=35939,274627$, the number of iterations for the ILU becomes up to $8$ times larger. For the unpreconditioned case, the maximum iteration count of $1000$ is reached for a significant number of linear system solves, limiting the accuracy of the resulting Newton steps.   

\paragraph{Drafting test with rectangular blade}

We follow the validation experiment in \cite{melanz2017comparison}, we set up a drafting test involving a rectangular blade of width $0.1$ moving through a container filled with spherical rigid particles of radius $0.1$ and friction coefficient $\mu = 0.25$. This test may be used to compute the force that the blade experiences as it moves through the granular flow, which eventually reaches a steady state. As in the sedimentation test above, we set up a perturbed lattice of $(2n+1)^3$ spherical particles inside the box, and perform experiments for $n=8,16$ for a total number of $M=4915,35939$ rigid bodies. The blade moves from one end of the box to the other with a prescribed sinusoidal velocity with period $4s$.    

\begin{figure}[H]
\centering
\includegraphics[width=0.7\textwidth]{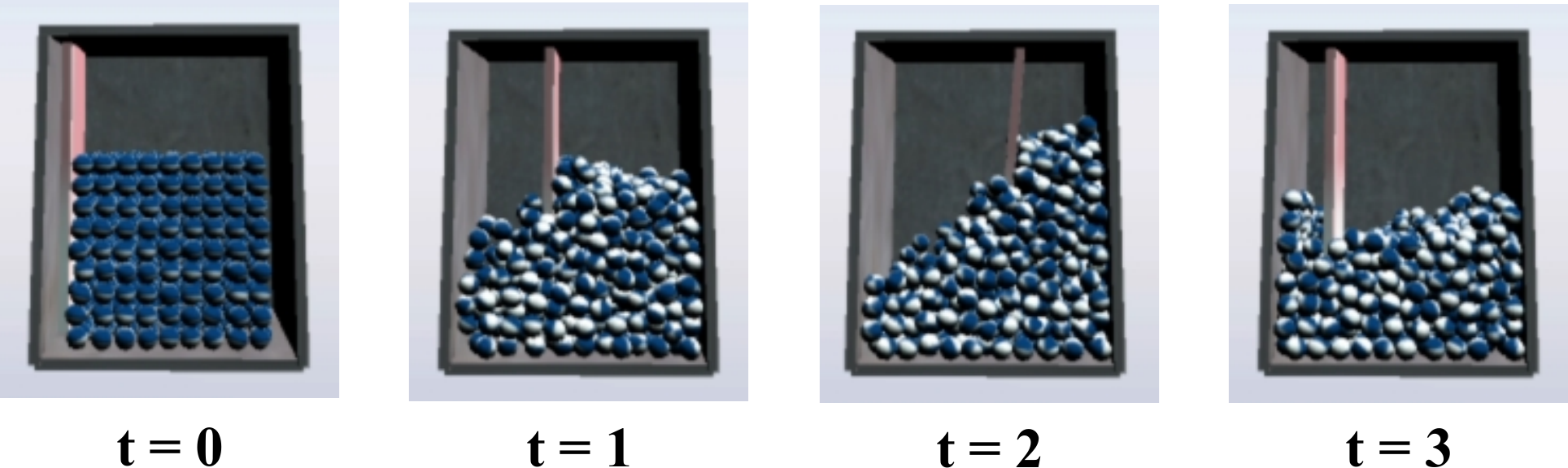}
\caption{\textbf{Blade drafting test:} snapshots of simulation for rectangular blade drafting test with $729$ spheres inside a closed box container. The blade moves on the $x$ direction with prescribed sinusoidal velocity.}
\label{fig:blade}
\end{figure}

\paragraph{Direct shear experiment} Finally, we include a direct shear test commonly used to measure the shear strength properties of granular soil, used in \cite{melanz2016experimental} to validate the DEM CCP approach for frictional contact. In this test, a soil sample is placed inside of a box and subjected to a normal load force (exterted by a cell weighing down on the material). The top half of the shear box is clamped while the lower is displaced in a controlled fashion from left to right, shearing the soil sample.This test can then be utilized to measure the shear stress and other relevant properties as a function of the shear displacement in the box. For our experiments, a perturbed box-shaped lattice with $(3n+1)(2n+1)^2$ spherical particles is set up inside the box, experiencing shear from the lower half of the box and a load from a ceiling press $10$ times denser than the granular material. We perform experiments for $n=8,16$, for a total number of $M=7228,53364$ rigid bodies. The movement of the lower half is again controlled prescribing a sinusoidal velocity. 

\begin{figure}[!htb]
\centering
\includegraphics[width=\textwidth]{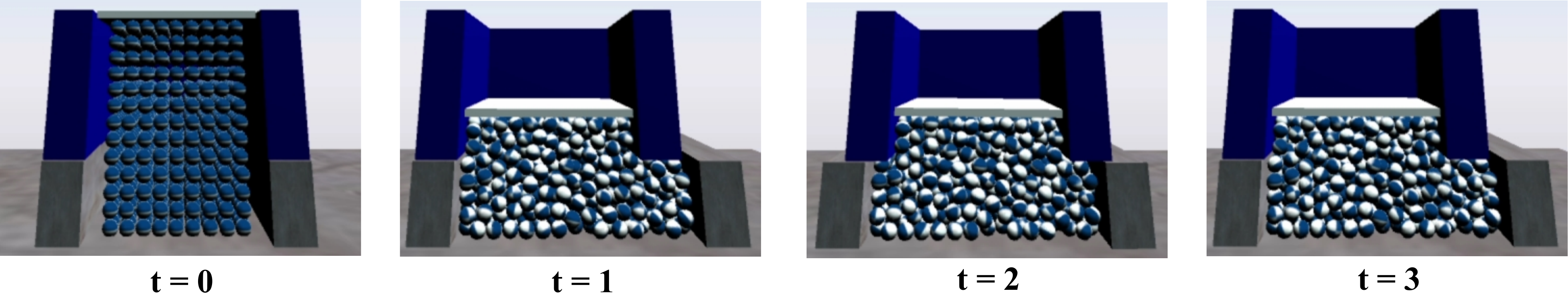}
\caption{\emph{Direct shear box experiment:} snapshots of simulation for direct shear test with $1053$ spheres inside a closed box container. The bottom half of the box is displaced in the $x$ direction with a prescribed sinusoidal velocity, and the granular fluid is loaded on the top by a rectangular press.}
\label{fig:shear}
\end{figure}

\begin{figure}[H]
\centering
\includegraphics[width=0.6\textwidth]{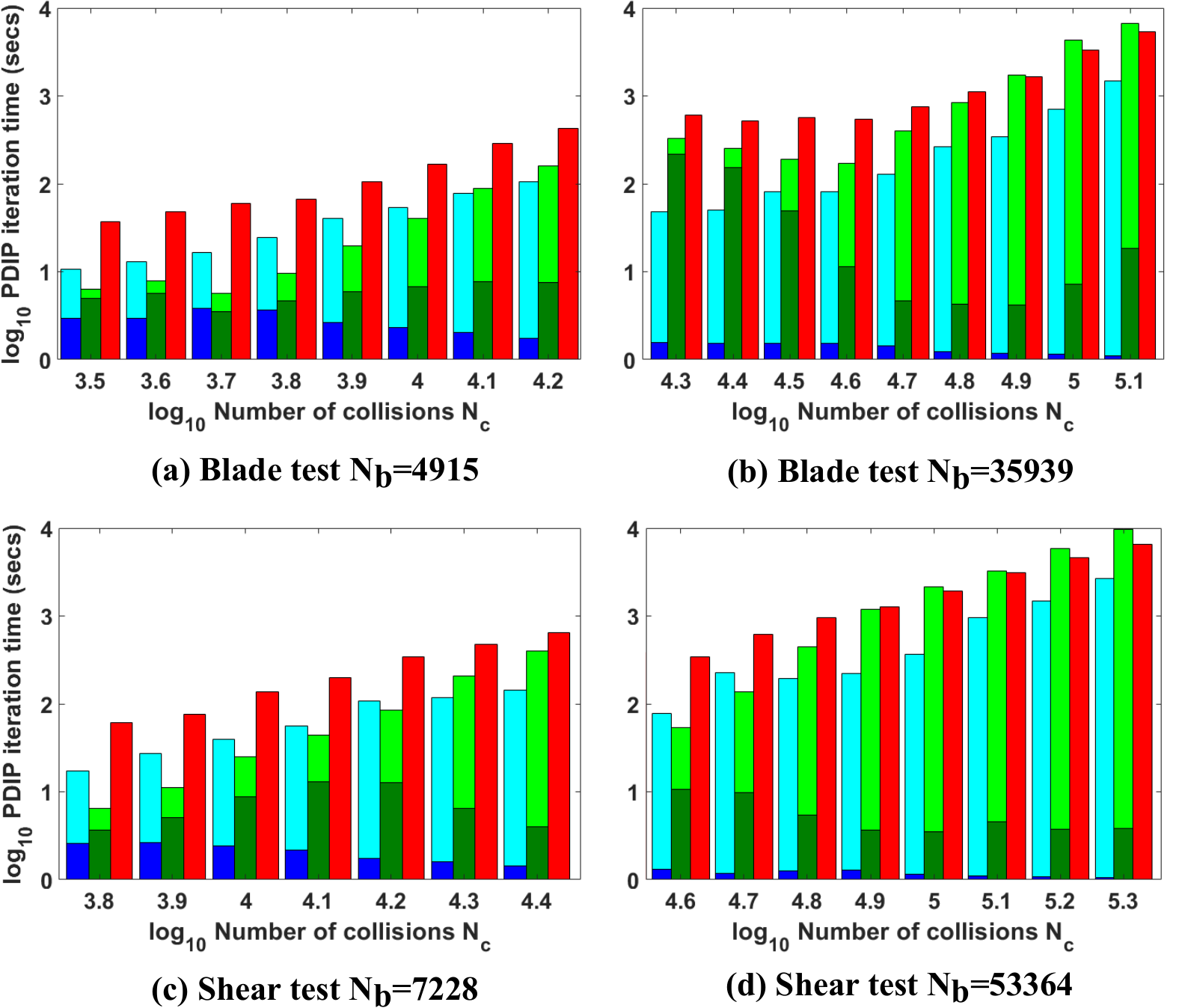}
\caption{\textbf{Log-log plot solver comparison for blade drafting (a,b) and direct shear (c,d) experiments:} for each experiment, we bin according to $log_{10} N_c$ and plot average PDIP iteration times; Unpreconditioned BICGSTAB in red, ILU-BICGSTAB in green and TT-BICGSTAB in blue. For preconditioned solvers, we display the proportion spent in precomputation in a darker shade, solve times in a lighter one.}
\label{fig:bladeshear_totstk}
\end{figure}

In \ref{fig:bladeshear_totstk}, we once again show a comparison of performance for each of the three solvers as it scales with number of collisions $N_c$. For all four experiments we observe essentially the same scaling for solution and precomputation times for the TT preconditioned solver, with TT precomputation times again staying roughly around 10 seconds per timestep. The overall speedups attained are slightly smaller (about 5-10x), due likely to the increase in problem complexity and the force and velocity magnitudes involved compared to the sedimentation case. However, it remains the case that the TT provides a significantly more robust and better acceleration to the linear system, with better scaling precomputation times and reduced iteration counts than the ILU sparse preconditioner. As is the case for the sedimentation tests, due to the limit of the maximum number of iterations for the unpreconditioned solve, average iteration counts and timings can be estimated to be about $5$ times higher than presented in these plots for full accuracy. 

\section{Conclusions}
\label{sec:conclusions}
In this work we have presented a robust and highly efficient acceleration technique for the solution of Newton step linear systems in second order methods based on approximate hierarchical compression and inversion in the Tensor Train format. In multiple experiments for common terramechanics phenomena modeled with frictional contact for dense, multiple rigid body systems, we have successfully applied a TT preconditioner to accelerate their solution, providing speed-ups of up to an order of magnitude against state-of-the-art sparse solvers for systems with $N_c \gtrsim 20000$. Across all our experiments, we observe that the Tensor Train preconditioner provides lower and more reliable iteration count reductions, as well as practically constant precomputation costs and storage requirements, which are improved significantly by our proposed TT factorization re-use techniques. 

As discussed in Section \ref{sec:intro}, the benefits of rapid, problem-independent convergence of second order optimization solvers are often negated by expensive large sparse matrix solves. The application of sparse and structured linear algebra techniques has been thus far limited by unfavorably scaling precomputation costs, excessive memory storage and communication requirements and the absence of efficient global factorization updates. We have demonstrated that the Tensor Train approach can successfully address these issues in DEM complementarity simulations of granular media. Based on its versatility and exploitation of a wide class of hierarchical low rank structure, we expect this to be true for a wide array of large scale optimization problems. We also anticipate that the low precomputation cost and storage requirements provided by the TT will be most impactful in high performance computing implementations; our ongoing work features a distributed memory implementation of the frictional contact CCP and the TT accelerated PDIP solver.    

\section{Acknowledgements}
We acknowledge support from the Automotive Research Center (ARC) in accordance with Cooperative Agreement W56HZV-14-2-0001 with U.S. Army Tank Automotive Research, Development and Engineering Center (TARDEC). Corona and Veerapaneni were also supported by the NSF under grant DMS-1454010.  This research was supported in part through computational resources and services provided by Advanced Research Computing Center at the University of Michigan, Ann Arbor. Corona would like to thank Daniel Negrut, Radu Serban, Luning Fang and Milad Rakhsha at the Simulation Based Engineering Laboratory (SBEL) at UW Madison for their assistance with this project and for many helpful discussions.  

\vspace{5mm}

\paragraph{\textbf{DISTRIBUTION STATEMENT A}}Approved for public release; distribution unlimited. OPSEC\#1118

\bibliography{corona}
\end{document}